\documentclass[]{amsart}

\usepackage{graphicx}
\usepackage{amsfonts}

\newtheorem{theorem}{Theorem}[section]
\newtheorem{lemma}[theorem]{Lemma}

\newtheorem{proposition}[theorem]{Proposition}

\theoremstyle{remark}
\newtheorem{remark}[theorem]{Remark}

\theoremstyle{definition}
\newtheorem{definition}[theorem]{Definition}

\newcommand{\gsigma}{\widetilde{\mathfrak{g}}_{\sigma}}
 \newcommand{\gsigmatau}{\widetilde{\mathfrak{g}}_{\sigma \tau}}

 \newcommand{\real}{{\mathbb R}}
 \newcommand{\rbold}{{\bf R}}
\newcommand{\lgsigma}{\Lambda G_{\sigma}}
\newcommand{\lgsigmatau}{\Lambda G_{\sigma \tau}}
 \newcommand{\lgsigmar}{\Lambda G_{\sigma}^{\bf R}}
\newcommand{\lgsigmataur}{\Lambda G_{\sigma \tau}^{\bf R}}

\newcommand{\DD}{{\mathbb D}_+}
\newcommand{\DC}{{\mathbb D}_-}
\newcommand{\chat}{\widehat{\mathbb C}} 
\newcommand{\betahat}{\hat{\beta}}

\newcommand{\rstar} {{\mathbb R}^*}
\newcommand{\cc}{{\mathbb C}}
\newcommand{\cstar} {{\mathbb C}^*}
\newcommand{\bbar}{\left[ \begin{array}}
\newcommand{\ebar}{\end{array} \right] }
\newcommand{\bdm}{\begin{displaymath}}
\newcommand{\edm}{\end{displaymath}}
\newcommand{\beq}{\begin{equation}}
\newcommand{\beqa}{\begin{eqnarray}}
\newcommand{\beqas}{\begin{eqnarray*}}
\newcommand{\eeq}{\end{equation}}
\newcommand{\eeqa}{\end{eqnarray}}
\newcommand{\eeqas}{\end{eqnarray*}}
\newcommand{\dd}{\textup{d}}
\newcommand{\hh}{\mathcal{H}}
\newcommand{\kk}{\mathcal{K}}

\begin{document}

\title[DPW method and an application to space forms]
{Generalized DPW method and an application to isometric immersions of space forms}   

\author{David Brander}
\address{Department of Mathematics\\ Faculty of Science\\Kobe University\\1-1, Rokkodai, Nada-ku, Kobe 657-8501\\ Japan}
\email{brander@math.kobe-u.ac.jp}

\author{Josef Dorfmeister}
\address{Tu Muenchen\\ Zentrum Mathematik (M8), Boltzmannstr. 3\\
  85748, Garching\\ Germany}
\email{dorfm@ma.tum.de}

\thanks{Research supported by DFG grant DO 776.}

\begin{abstract}
Let $G$ be a complex Lie group and $\Lambda G$ denote the group of 
maps from the unit circle
${\mathbb S}^1$ into  $G$, of a suitable class.
A differentiable map $F$ from a manifold $M$ into $\Lambda G$,
is said to be of \emph{connection order $(_a^b)$} if the  
Fourier expansion
in the loop parameter $\lambda$ of the ${\mathbb S}^1$-family of Maurer-Cartan forms for
$F$, namely $F_\lambda^{-1} \dd F_\lambda$, is of the form
$\sum_{i=a}^b \alpha_i \lambda^i$. Most integrable systems in
geometry are associated to such a map.
Roughly speaking, 
the DPW method used a Birkhoff type splitting to reduce a harmonic
map into a symmetric space, which can be represented by a certain order $(_{-1}^1)$ map,
 into a pair of simpler
 maps  of order $(_{-1}^{-1})$
and $(_1^1)$ respectively.
Conversely, one could construct such a harmonic
map from any pair of  $(_{-1}^{-1})$ and  $(_1^1)$  maps. 
This allowed a Weierstrass type description of harmonic
maps into symmetric spaces. We extend this method to show that, for a 
large class of loop groups,
a connection order $(_a^b)$ map, for $a<0<b$, splits uniquely into a pair of 
$(_a^{-1})$ and $(_1^b)$ maps.  
  As an application, we show that constant non-zero sectional
curvature submanifolds with flat normal bundle of a sphere or hyperbolic space 
split into pairs of flat submanifolds, reducing the  problem
(at least locally)
to the flat case.  To extend the DPW method sufficiently
 to handle this problem
requires a more general Iwasawa type splitting of the
 loop group, which we prove
always holds at least locally.
\end{abstract}


\maketitle

\section{Introduction} \label{intro}
Over the past two or three decades,
various techniques based on the Birkhoff and Iwasawa loop group
factorizations
have been used successfully to study integrable systems
in geometry. This paper may be viewed as an attempt to clarify generally the
class of problems to which they apply, in addition to considering a new
application. 

\subsection{Limited connection order maps into loop groups} \label{sublcom}
Let $G$ be a complex Lie group, with 
Lie algebra $\mathfrak{g}$, and $\Lambda G$ the group of maps from the 
unit circle ${\mathbb S}^1$ into $G$, with a topology that makes $\Lambda G$
a Banach Lie group which has some additional properties needed for the purposes of this paper.
 For a matrix group $G$, the Wiener topology has all these
properties, and includes all applications the authors are aware of.

  Let $M$ be a smooth manifold and $F$ a smooth function from
$M$ into $\Lambda G$. Denote the group of such maps by $\Lambda G(M)$.
 $F$ is often regarded as a family of smooth maps $F_{\lambda}$ from
$M$ into $G$, parameterized by $\lambda \in {\mathbb S}^1$.   
In this sense ``the Maurer-Cartan form for $F$" is
the family $A_\lambda := F_\lambda^{-1} \dd F_\lambda$
 of Maurer-Cartan forms for
$F_{\lambda}$.  The integrability condition for $A_\lambda$
 is the Maurer-Cartan equation
\beq \label{mceqn}
\dd A_\lambda + A_\lambda \wedge A_\lambda = 0,
\eeq
and this is satisfied for all $\lambda$ in ${\mathbb S}^1$.  Conversely, consider a family
$A_{\lambda}$ of 1-forms on $M$, with values in $\mathfrak{g}$, 
such that 
the map from $M$ to the Banach Lie algebra $\textup{Lie} \Lambda G =: \Lambda \mathfrak{g}$, 
associating to $p \in M$ 
the ${\mathbb S}^1-$family  $A_\lambda (p)$, is smooth
 and
 which satisfies, in addition, (\ref{mceqn}) for all $\lambda$. Then, on
any simply connected domain in $M$, there exists a family $F_{\lambda} \in \Lambda G(M)$
whose Maurer-Cartan form is $A_{\lambda}$.

If we expand $A_\lambda$ as a Fourier series in $\lambda$,
\bdm
A_{\lambda} = \sum_i a_i \lambda ^{i},
\edm
where each $a_i$ is a $\mathfrak{g}$-valued 1-form on $M$, then the equation 
(\ref{mceqn}) is equivalent to the system of equations
\beq  \label{intsystem}
\dd a_k + \sum_{i+j = k} a_i \wedge a_j = 0.
\eeq
By assuming that $F$ is a map into a particular subgroup $\hh$ of $\Lambda G$, for
 some
particular group $G$, and that only a few $a_i$ are non-zero, it has
been found that certain maps of interest in geometry - such as harmonic maps 
into symmetric spaces, \cite{hitchin1986}, \cite{pinkallsterling},
 \cite{fpps}, \cite{uhlenbeck1989}, \cite{bfpp}, \cite{burstall1995},
 isometric immersions of 
space forms into space forms \cite{feruspedit1996}, \cite{feruspedit1996II},
\cite{lms},  Hamiltonian stationary Lagrangian surfaces in
${\mathbb C}^2$ \cite{heleinromon} -
 can be characterized by such an $F$, where the 
the special
equations defining  the particular type of map are precisely 
(\ref{intsystem}).  A recent survey which catalogues most of the known
examples is given by Terng \cite{terng2002}.
More specific information concerning the construction of surface classes can be found in 
\cite{Dorfts}.

We will say that an element $F \in \Lambda G(M)$ is of 
\emph{connection order $(_a^b)$}, where $a$ and $b$ are
extended integers and $a\leq b$, if the Maurer-Cartan form of $F$ is a Laurent series of
top and bottom degree $b$ and $a$ respectively,
\bdm
F^{-1} \dd F = \sum_a^b a_i \lambda^{i},
\edm
and denote the subset of such maps by $\Lambda G(M)_a^b$. These are the 
basic objects of study in this article.

Let $\Lambda^+ G$ and $\Lambda^- G$ denote the subgroups of $\Lambda G$
consisting of loops which extend to holomorphic maps into $G$ from the 
unit disc and its complement in the Riemann sphere respectively.
Involutions of the loop group $\Lambda G$ can be divided into two basic types, those
of the \emph{first} and those of the \emph{second kind}. These are essentially characterized
by the property that they map $\Lambda^\pm G$ to $\Lambda^\pm G$ for the first kind,
and $\Lambda^\pm G$ to $\Lambda^\mp G$ for the second kind.
An important subclass of limited connection order maps are the
\emph{connection order $(_{-b}^b)$ maps with involution of the second kind $\tau$}, or
$(_{-b}^b)_\tau$ maps, where $b \geq 0$.
These are  maps into 
a subgroup $\Lambda G_\tau$, 
the fixed point set of an
 involution of the second kind $\tau$.  
In most applications that the authors
are aware of, one seeks to construct maps into a real form $G^\rbold$ of
the complex Lie group $G$, and this leads to consideration of loops
(actually maps from from ${\mathbb C}^*$ into $G$) which are real along
either a line through the origin or along ${\mathbb S}^1$.  The ${\mathbb S}^1$ reality
condition is of the second kind, and so leads to type $(_{-b}^b)_\tau$ maps;
a different example of a $(_{-b}^b)_\tau$ map will be discussed below, in 
Section \ref{spaceforms}.
 
For reasons which will become clear, we distinguish between two types of
loop group problems: those with and those without an involution of the second kind.
\subsection{Examples without an involution of the second kind} \label{sub1}
The $(b,j)$'th flow of the $G$-hierarchy \cite{terng2002} originating
in the well known ZS-AKNS  construction  \cite{zs}, \cite{akns}, 
is associated to an order $(_0^j)$ map. The $-1$-flow and the 
hyperbolic $U$-system \cite{terng2002} are $(_{-1}^1)$  maps.
Some known examples of geometrical problems of such type are 
pseudospherical surfaces in $\real ^3$ \cite{todathesis}, \cite{toda2002},
\cite{toda2005}, and
affine spheres \cite{dorfeitner}, both of order $(_{-1}^1)$, and
curved flats \cite{feruspedit1996II},  which are of order $(_0^1)$.
Particular examples of curved flats are flat isometric immersions into
a sphere or hyperbolic space \cite{lms} and isothermic surfaces \cite{bhpp}.

\subsection{Examples with an involution of the second kind} \label{sub2}
Harmonic maps into
symmetric spaces are of connection order $(_{-1}^1)_\tau$, and, more
generally, the $m$'th elliptic $(G, \tau)$ and $(G, \tau, \sigma)$ 
systems defined by Terng \cite{terng2002} are of order $(_{-m}^m)_\tau$.
Some concrete examples are Hamiltonian stationary
Lagrangian surfaces in ${\mathbb C}^2$, which  are of  order $(_{-2}^2)_\tau$ \cite{heleinromon}, 
and isometric immersions with flat normal 
bundle of \emph{non-flat} space forms into a sphere or hyperbolic space, which are of  order
 $(_{-1}^1)_\tau$ \cite{feruspedit1996} (see also Section \ref{spaceforms} below).

\subsection{The Birkhoff factorization theorem} \label{introbirkhoff}
The basis of the DPW method is  the Birkhoff factorization theorem for $\Lambda G$
(Theorem \ref{birkhoff}), which  states that: each element of $\Lambda G$  can be written in the form
$g = g_- w g_+$, where $g_{\pm} \in \Lambda^{\pm} G$  and $w$ is in a specific set (which is in 
most of the cases of interest to us the Weyl group of the Lie algebra $\Lambda \mathfrak{g}$).
Moreover, the  \emph{left big cell},  
$\mathcal{LB} \Lambda G  = \Lambda G^- \cdot \Lambda G ^+$ is an 
open and dense subset of $(\Lambda G)_I$, where $(X)_I$ denotes the identity component
of a group $X$. By definition, $\mathcal{LB} \Lambda G $ 
consists exactly of those elements of $g \in \Lambda G$, which can be written in the form 
$g = g_- g_+$, where $g_{\pm} \in \Lambda G^{\pm}$. 
We note that such a decomposition can be normalized in various ways.
Usually we use $g_+ ( 0) = I$ or $g_- (\infty ) = I$ to make the decomposition unique.
In this case, $g_+$ and $g_-$ are analytic functions of $g$. 
There is obviously an analogous statement with $\pm$ interchanged. 
The corresponding {\it right big cell} will be denoted $\mathcal{RB} \Lambda G $.
The intersection  
$\mathcal{B} \Lambda G  =    \mathcal{LB} \Lambda G  \cap \mathcal{RB} \Lambda G $
will be of particular importance. It will be called the {\it (central) big cell}.

For certain purposes it will also turn out to be useful to consider  groups more general
than the loop groups discussed above. In particular,
if such a group is a  subgroup $\hh$ of $\Lambda G$, we define  subgroups 
 $\hh^{\pm} = \hh \cap \Lambda ^\pm G$.

We will call a subgroup $\hh$ of $\Lambda G$ {\it Birkhoff decomposable} if  
the (central) big cell
$\mathcal{B} \hh := ((\hh^+)_I \cdot (\hh^-)_I) \, \cap \,  ((\hh^-)_I \cdot (\hh^+)_I)$ is open and dense in $(\hh)_I$. Well known examples of 
Birkhoff decomposable groups are  
the $\sigma-$twisted subgroups ${\Lambda G}_{\sigma}$, which denotes a fixed point 
subgroup of $\Lambda G$, with respect to an involution of the first kind. 
But we will also present, for example in Section \ref{birkhoffsubsect}, 
Birkhoff decomposable groups different from those just mentioned.

In most cases, as considered in this paper, the loop groups which have arisen in integrable systems are
either Birkhoff decomposable or a fixed point subgroup $\hh_\tau$, where $\hh$
is Birkhoff decomposable and $\tau$ is an involution of $\hh$ of the second kind.

  
\subsection{The generalized DPW method}  \label{sub3}
The theme of this paper is to use Birkhoff decompositions to break down
connection order $(_a^b)$ maps (where $a<0<b$) into pairs of $(_a^{-1})$
and $(_1^b)$ maps.  This generalizes the DPW method \cite{dorfmeisterpeditwu},
where the idea was applied to harmonic maps, an order $(_{-1}^1)_\tau$ example.
In fact the pair of $(_1^1)$ and $(_{-1}^{-1})$ maps obtained in that
example were related by the involution $\tau$, so the solution could be described by a
single order $(_1^1)$ map, whose Maurer-Cartan form  turned out to be a meromorphic Lie algebra valued function of one complex variable.

 We first prove (Propositions \ref{prop1} and \ref{prop2}) that if $a<0<b$, and $\hh$ is any Birkhoff decomposable subgroup,
then connection order $(_a^b)$ maps into $\hh$ are
in one to one correspondence with pairs of connection order $(_a^{-1})$
 and $(_{1}^b)$ maps. More accurately, the converse part, Proposition \ref{prop2},
 allows one to construct an order $(_a^b)$ map from a pair of $(_a^x)$ and $(_y^b)$
 maps, where $x \geq a$ and $y \leq b$ are arbitrary extended integers.
Propositions \ref{prop1} and \ref{prop2} are a straightforward generalization of 
 similar results in \cite{todathesis},  where 
they were
proved for a loop group associated to 
 pseudospherical surfaces in $\real ^3$. In that case, the resulting
pair of $(_{-1}^{-1})$ and $(_1^1)$ maps were each a Lie algebra valued function of
one variable only,
 in asymptotic coordinates,
 which enabled a simplified ``Weierstrass-type''
characterization of pseudospherical surfaces. A similar outcome was obtained
in \cite{dorfeitner} for affine spheres. 
Essentially the same idea had been used by Krichever in \cite{krichever},
before the introduction of the loop group point of view.

As noted above, many of the known applications are of order $(_{-b}^b)_\tau$,
that is they map into a subgroup $\hh_{ \tau}$, where $\hh$ is Birkhoff decomposable and where $\tau$ is
an involution of $\hh$ of the second kind.
It is easy to see that $\hh_{\tau}$ is not Birkhoff decomposable: this makes the converse
part of the splitting tricky. We will show that the problem can be overcome when $\hh$ is 
\emph{$\tau$-Iwasawa decomposable}, that is, if 
a neighbourhood of the identity of $\hh$ has a unique
factorization 
\bdm
F_\tau  g_+,
\edm
where $F_\tau$ is an element of $\hh_{ \tau}$ and $g_+$ is in $\hh^+$.
For such a situation, we 
prove (Proposition \ref{prop3}) 
that, in fact, there is a correspondence between connection order $(_{-b}^b)_\tau$
 maps and single order $(_1^b)$ maps. 
 
It turns out  (Theorem \ref{iwasawathm}) that
$\hh$ is always $\tau$-Iwasawa
decomposable, for any Birkhoff decomposable subgroup $\hh$,
and we can make more precise statements in the case that the subgroup of constant loops, $\hh^\circ$, 
is \emph{compact}.  

We note here (see Remark \ref{doubleloopgroup}),
that, by considering the double loop group of \cite{dowu},
  one can recast these results
so that one only need consider involutions of the first kind.

In Section \ref{dressing} we include a short discussion on how
dressing defines a natural local action by $\hh ^-$ 
on $\hh(M)_1^b$, and similarly an action by $\hh ^- \times \hh ^+$
on $\hh(M)_a^b$.

In the last sections, we apply the splitting results to the cases
of isometric immersions with flat normal bundle of space forms
$M_c \to \widetilde{M}_{\tilde{c}}$, where $c$ and $\tilde{c}$ 
are the respective constant sectional curvatures.
This problem had earlier been studied in the loop group framework,
in \cite{feruspedit1996}, for non-flat immersions,
and in  \cite{lms} and \cite{terng2002} for the flat case.

We show that, given a choice of base point for the manifold,
 the DPW splitting produces 
(at  least at a local level) 
a 1-1-relation between all 
non-flat immersions and all flat immersions, in accordance with Table \ref{cortable}.
\begin{table}[here]
\caption{Local correspondence between non-flat and flat immersions, $M_c \to \widetilde{M}_{\tilde{c}}$}
\label{cortable}
  \begin{tabular}{|c|c||c||c|c|}  \hline  
$c \in (-\infty,0)$ & $\tilde{c} = 1$ & $\leftrightarrow$ & $c = 0$ & $\tilde{c} = 1$ \\ \hline
$c \in (0,1)$ & $\tilde{c} = 1$ & $\leftrightarrow$ & $c = 0$ & $\tilde{c} = 1$ \\
\hline
$c \in (1,\infty)$ & $\tilde{c} = 1$ & $\leftrightarrow$ & $c = 0$ & $\tilde{c} = -1$ \\ \hline
 $ c \in (-\infty,-1)$ & $\tilde{c} = -1$ & $\leftrightarrow$ & $c = 0$ & $\tilde{c} = 1$ \\ \hline
$c \in (-1,0)$ & $\tilde{c} = -1$ & $\leftrightarrow$ & $c = 0$ & $\tilde{c} = -1$ \\
 \hline
  $c \in (0,\infty)$ & $\tilde{c} = -1$ & $\leftrightarrow$ & $c = 0$ & $\tilde{c} = -1$ \\
     \hline
  \end{tabular} 
\end{table}

As another application of the generalized DPW method
described here, the technique is used, in the article
\cite{brander2}, to relate isometric immersions of space forms
to pluriharmonic maps, and that particular setup is generalized 
further to constant curvature Lagrangian submanifolds of space forms
in \cite{reflective}.  

\subsection{Further directions}
Certain special connection order $(_a^b)$ maps of interest 
are not treated here.
For example, so-called ``finite type'' solutions,
a description of which can be found in \cite{burstallpedit}.
It is not difficult to show that for  a finite type connection order
$(_a^b)$ map, for $a<0<b$, 
the corresponding pair of $(_a^{-1})$ and $(_1^b)$ maps
are also of finite type.
Another interesting class of maps which may merit 
investigation are solutions with ``finite uniton number", studied by
Uhlenbeck \cite{uhlenbeck1989} - see also Guest and Ohnita \cite{guestohnita}.

\section{Birkhoff and $\tau$-Iwasawa decompositions for subgroups of $\Lambda G$}
\label{birkhoffiwasawa}
In the next few sections, we discuss some general theory of
(loop)  group decompositions and their applications to 
connection order $(_a^b)$ maps. Throughout this discussion,
$\Lambda G$ will generally denote a
loop group, and $\hh$ some Birkhoff decomposable subgroup of $\Lambda G$.

\subsection{Definitions and terminology}
 For convenience, we gather together here some definitions.
\subsubsection{The loop group} 
All groups of interest to us can be realized as matrix groups, so we proceed as
follows: let  $G$ be a connected complex
semisimple matrix Lie group, with some matrix norm $\| \cdot \|$,
satisfying $\|AB \| \leq \|A \| \| B \| $ and   $\| I \| = 1$,
and with Lie algebra $\mathfrak{g}$. Define
the loop group
\bdm
\Lambda G := \{ g : {\mathbb S}^1 \to G ~ |~ g = \sum_{i=-\infty}^\infty g_i \lambda^i, ~~ ~ \sum \|g_i\| < \infty \}. 
\edm

Equipped with  the  topology induced by the norm 
$\|g\| := \sum \|g_i\|$,  $\Lambda G$ is a Banach Lie group whose Lie algebra,
$\Lambda \mathfrak{g}$,  consists of Laurent series in $\lambda$ with 
coefficients in the Lie algebra of $G$, and the corresponding convergence condition.
For more details see \cite{balandorf} or \cite{pressleysegal}.

\subsubsection{Subgroups}
All subgroups of $\Lambda G$ are assumed to be Banach Lie subgroups in the sense of 
\cite{bourbaki}.
In particular, subgroups are
closed submanifolds and carry the induced
topology. 
Let $\chat$ be the extended  plane,  $\DD$ the open unit disc, and
$\DC$ the complement of its closure in $\chat$.
For any group $G$ as above, define the following subgroups of $\Lambda G$:
\bdm
\Lambda ^{\pm}G := \{ g \in \Lambda G~ | ~g \textup{ has a holomorphic extension } \hat{g}: {\mathbb D}_{\pm} \to G  \}. 
\edm
These have the subgroups $\Lambda ^{+}_1G$ and $\Lambda ^{-}_1G$ whose elements have
the respective normalizations  $g (0) = I$ and $g (\infty) = I$.

For a subgroup $\hh$ of $\Lambda G$, we define $\hh^{\pm}$ and $\hh^{\pm}_1$ to be the
intersection of $\hh$ with the notationally corresponding subgroups of $\Lambda G$.
In all concrete examples we consider, all these groups  will be 
Banach lie groups and subgroups respectively.

If $\sigma$ and $\rho$ are automorphisms of $\hh$, then $\hh_{\sigma \rho}$ denotes the
subgroup of $\hh$ consisting of elements fixed by both automorphisms, and similar notation
is used for any number of automorphisms. In particular, $\Lambda G _{\sigma}$ means
$(\Lambda G)_\sigma$.

\subsection{Birkhoff factorization} \label{birkhoffsubsect}
The Birkhoff factorization theorem, as stated here, 
is proved, for the case of $\mathcal{C}^\infty$ loops
 in \cite{pressleysegal} for $GL(n, {\mathbb C})$ as well as for
 any complex semisimple Lie group. In \cite{balandorf}, a similar result is shown
to hold, working with the Wiener topology,
 for the complexification of any connected Lie group which admits
a faithful finite dimensional continuous representation.

\begin{theorem} \label{birkhoff}
\begin{enumerate}
\item
Any loop $g \in \Lambda G$ has a 
\emph{right Birkhoff factorization}
\beq \label{rightfact}
g = g _+ D g _-,
\eeq
where $g _+ \in \Lambda ^+ G$, 
$g _- \in \Lambda ^- G$ and $D$ is a homomorphism from ${\mathbb S}^1$
into a maximal torus of $G$.
Similarly a  \emph{left Birkhoff factorization}, 
\beq \label{leftfact}
g = g _- D g _+,
\eeq
exists, where $g _- \in  \Lambda ^- G$ and
 $g _+ \in  \Lambda ^+G$.
\item
The \emph{left big cell}, 
$\mathcal{LB} \Lambda G$, defined to be the subset on which
the left factorizations have $D = I$ as the center term,
is open and dense in $(\Lambda  G)_I$. In fact it is 
given as the complement of
the zero set of a non-constant 
holomorphic section $f_L$, of a holomorphic line
bundle over $(\Lambda G)_I$ (see \cite{dorfmeisterpeditwu}).
The \emph{right big cell} has the analogous definition and properties.
\item
 The maps  
$\Lambda ^+_1 G \times  \Lambda ^- G \to \mathcal{RB} \Lambda G$, and
$\Lambda ^-_1 G \times  \Lambda ^+ G
\to \mathcal{LB} \Lambda G$ are analytic diffeomorphisms as are the maps 
$\Lambda ^+ G \times  \Lambda ^-_1 G \to \mathcal{RB} \Lambda G$ and
$\Lambda ^- G \times  \Lambda ^+_1 G \to \mathcal{LB} \Lambda G$.
In particular, on the big cells the factorization is unique, if one requires 
$g_+ \in \Lambda^+_1 G$ or $g_- \in \Lambda^-_1 G$.
\end{enumerate}
\end{theorem}

 The intersection,
 $\mathcal{B} \Lambda G := \mathcal{LB}\Lambda G\cap  \mathcal{RB}\Lambda G$, we will call the
 \emph{(central) big cell}.

More generally, a Banach Lie group $\hh$ together with connected Banach Lie subgroups 
subgroups $\hh_1^{\pm}$ and 
${\hh}^\circ$
will be called {\it Birkhoff decomposable}, if 
$\mathcal{RB} \hh := \hh_1^+ \cdot {\hh}^{\circ} \cdot \hh_1^-$ and
$\mathcal{LB} \hh := \hh_1^- \cdot {\hh}^{\circ} \cdot \hh_1^+$ are 
open and dense in $(\hh)_I$ and the maps
$ \hh_1^+ \times {\hh}^{\circ} \times \hh_1^- \rightarrow \mathcal{RB} \hh$ and 
$ \hh_1^- \times {\hh}^{\circ} \times \hh_1^+ \rightarrow \mathcal{LB} \hh$ are 
analytic diffeomorphisms.
In this case every element $h \in \mathcal{RB} \hh$ has a unique decomposition 
$ h = h_1^+ \cdot {h}^{\circ} \cdot h_1^- $ with $h_1^{\pm} \in \hh_1^{\pm}$ and 
${h}^{\circ} \in
{\hh}^{\circ}$. The analogous statement for elements 
$h \in \mathcal{LB} \hh$ also holds.
Setting $\hh^+ = \hh_1^+ \cdot {\hh}^{\circ}$ and 
$\hh^- = \hh_1^- \cdot {\hh}^{\circ}$ we obtain the 
properties listed in the last part of Theorem \ref{birkhoff}. 

If $\hh$ is a subgroup of $\Lambda G$, then usually the subgroup
$\hh^\circ$ referred to above is just the constant subgroup
$\hh \cap G$; in the general case,
we allow other possibilities a priori.  In this article,
however, all examples of Birkhoff decomposable subgroups considered 
will be related in one way or another to
subgroups of loop groups. 

Let $\hh \subset \Lambda G$ be a Banach Lie subgroup.
\begin{definition}
We will call $\hh$ a Birkhoff decomposable subgroup of $\Lambda G$ if 
the identity components of
$\hh^+_1 := \hh \cap \Lambda G^+_1$ , $\hh^-_1 := \hh \cap \Lambda G^-_1$ and
${\hh}^{\circ} := \hh \cap \Lambda^+ G \cap \Lambda^- G$ are closed Banach Lie
subgroups and make $\hh$ 
into a Birkhoff decomposable group.
\end{definition}

{\it Example 1 } The groups $\Lambda G$ are Birkhoff decomposable groups relative to
$\Lambda G^+_1$, $\Lambda G_1^-$ and ${\Lambda G}^{\circ} = \Lambda G^+ \cap \Lambda G^- = G$.

{\it Example 2 } Let ($\hh$, $\hh_1^+$, ${\hh}^{\circ}$, $\hh_1^-$) and 
($\kk$, $\kk_1^+$, ${\kk}^{\circ}$, $\kk_1^-$) be Birkhoff decomposable groups.
Then the product group, together with the obvious subgroups is also a Birkhoff decomposable group.

{\it Example 3 } Let ($\hh$, $\hh_1^+$, ${\hh}^{\circ}$, $\hh_1^-$) be a 
Birkhoff decomposable 
subgroup of some loop group $\Lambda G$.
\begin{theorem} \label{productthm}
Set
\beqas
\kk := \hh \times \hh, &\hspace{1cm}& {\kk}^{\circ} = {\hh}^{\circ} \times {\hh}^{\circ},\\
\kk_1^+ = \hh_1^+ \times \hh_1^- &\hspace{1cm}&
\kk_1^- = \lbrace (x,x);~ x \in \hh,~ x|_{\lambda = 1} = I \rbrace.
\eeqas
Then $\kk$ with the corresponding subgroups is a Birkhoff decomposable group.
\end{theorem}
\begin{proof}
 We note first that $\kk$ is a Banach Lie group and the remaining three 
groups are Banach Lie subgroups of $\kk$. For $\kk_1^-$ one uses the fact that for  
loop groups with the Wiener topology, the evaluation map is an 
analytic homomorphism of Lie groups, such that its kernel is a Banach Lie subgroup 
(actually a normal subgroup, of course).
The decisive object now is the map: $\Psi: \kk_1^- \times {\kk}^{\circ} \times \kk_1^+
 \to \kk$ given by
 $((h,h),(a,b), (g_+,g_-)) \mapsto (hag_+,hbg_-)$.
It is straightforward to prove that it is analytic and injective. Considering 
the inverse map we start from
$(p,q) = (hag_+,hbg_-)$. Then $p^{-1}q = g_+^{-1}a^{-1}bg_-$. 
Therefore, using that $\hh$ is a Birkhoff decomposable group,
$g_+$ and $g_-$ are analytic in $p$ and $q$ as well as $a^{-1}b$. But then the equation 
$ pg_+^{-1} = ha$ shows, by evaluating at $\lambda =1$, that $a \in G$ is analytic in $p$ and $q$.
Similarly one sees that $b$ is analytic in $p$ and $q$. Thus, finally, also $h$ is analytic in $p$ and $q$.

It remains to check that the image $\textup{Im} \, \Psi$ of $\Psi$ is open and dense in $(\kk)_I$.
This can be verified by observing that $(p,q)$ lies in $\textup{Im} \, \Psi$  if and only if we
can write $p = hf_+$ and $q = hf_-$, for some $h \in \hh$ and $f_\pm \in \hh^\pm$. 
This is equivalent to solving $pf_+^{-1} = q f_-^{-1}$, or 
$q^{-1}p = f_-^{-1}f_+$, which has a solution if and only if $q^{-1}p \in \mathcal{LB}\hh$.
Therefore, the image of $\Psi$
 is $\phi ^{-1} (\mathcal{LB}\hh)$, where $\phi(p,q) := q^{-1}p$ is a real analytic map,
and $(p,q) \in \textup{Im} \, \Psi$ if and only if $\phi^* f_L (p,q) \neq 0$,
where $f_L$ is the holomorphic section given in Theorem \ref{birkhoff}.
But $\phi^* f_L$ is real analytic, and so the complement of its zero set
is either empty or open and dense. It cannot be empty, as it contains the identity.
\qed
\end{proof} 
We will give more examples of Birkhoff decomposable groups in the following sections.

\subsection{Automorphisms and Birkhoff decomposable subgroups}
Let $\sigma: \Lambda G \to \Lambda G$ be a real analytic   automorphism, with  
fixed point subgroup $\lgsigma$.
Clearly, we can Birkhoff factorize  elements in $\Lambda G_{\sigma}$ with factors in 
the larger group $\Lambda G$.  We would like the factors to be in $\Lambda G_\sigma$
themselves, more strongly even,  the subgroup to be 
 Birkhoff decomposable.
More generally, if $\hh$ is any Birkhoff decomposable subgroup we ask the same
question for fixed point subgroups $\hh_\sigma$ of $\hh$.

Finite order automorphisms of $\Lambda G$ have been classified \cite{levstein}, \cite{kacpeterson},
\cite{bauschrousseau}. See also \cite{heintze}.
After passing to an isomorphic loop group, one can assume, for a complex linear automorphism
$\phi$, of finite order $n$, that it  can be expressed as 
\beq \label{autocomp}
(\phi x)(\lambda) = \phi_0 (x(e^{\frac{2 k \pi i}{n}} \lambda^j)),
\eeq
where  $\phi_0$ is a finite order automorphism of $G$, $k$ is an integer,  and  $j = \pm 1$.
If $j$ is positive, the automorphism is said to be \emph{of the first kind}.
If negative, \emph{of the second kind}. In the latter case one can even assume that the factor 
$e^{\frac{2k \pi i}{n}}$ is not present.
However, for simplicity of presentation we will not use this possibility in general.

Complex \emph{anti}-linear automorphisms  of finite order can be classified
similarly. In particular, if $\rho$ is an anti-linear involution
 then one can assume it is of the form
\beq \label{rhopm} 
(\rho x)(\lambda) := \rho_0 (x(\varepsilon (\bar{\lambda})^{j})),
\hspace{1cm} \varepsilon ,~j \in \{1,-1\},
\eeq
for some anti-linear involution $\rho_0$ of $G$. These are also said to be of the
first and second kind if $j$ is positive or negative respectively.
As above, for involutions of the second kind the factor $\varepsilon$ can be removed, 
but it is convenient not to do this.

\begin{remark} In fact there are applications where one is also interested in automorphisms
which are not in standard form, such as the $n \times n$ KdV reality condition used in
\cite{ternguhlenbeckkdv}.  In this article, however, all automorphisms are assumed
to be of the forms given by (\ref{autocomp}) and (\ref{rhopm}).
We note that all the automorphisms of finite order used in this paper are analytic and their fixed point subgroups are Banach Lie subgroups.
\end{remark}

A \emph{real form} of $\Lambda G$ is defined to be the fixed point
set of an anti-linear involution.

At the level of Lie algebras, one also has expressions of the form (\ref{autocomp}) and
(\ref{rhopm}), replacing $\phi_0$ and $\rho_0$ with the corresponding Lie algebra
automorphisms.

Note that an automorphism of the first kind takes $\Lambda ^{\pm}G \to \Lambda ^{\pm}G$,
while one of the second kind  takes $\Lambda ^{\pm}G \to \Lambda ^{\mp}G$.

\begin{proposition} \label{generaltwistedfact}
Let $\hh$ be a Birkhoff decomposable subgroup of $\Lambda G$,  and
suppose $\phi$ is a finite order
 automorphism of $\Lambda G$ which restricts to an automorphism of $\hh$. 
Then:
\begin{enumerate}
\item
If $\phi$ is of the first kind, then $\hh_\phi$ is Birkhoff decomposable.
\item
If $\phi$ is of the second kind, and $\hh_\phi$ does not consist entirely of constant loops,
then $\hh_\phi$ is not Birkhoff decomposable.
\end{enumerate}
\end{proposition}
\begin{proof}
\begin{enumerate}
\item
We give the proof for the right factorization here, the left
being analogous. Suppose $x \in \mathcal{RB}\hh \cap \hh_\phi$.
Since $x$ is fixed by $\phi$, we have  $x = \phi x$.
Since $x \in  \mathcal{RB} \hh$, $x$ has a unique decomposition
$ x = x_+ {x}^{\circ} x_-$ with $x_{\pm} \in \hh_1^{\pm}$ and 
${x}^{\circ} \in \hh^{\circ}$.
First we show that these factors are already contained in $\hh_\phi$ and 
that the required density and analyticity properties hold.
First we note
\beq \label{eqn123}
x_+  {x}^{\circ} x_-  = 
\phi x_+ \phi {x}^{\circ} \phi x_-.
\eeq
Now $\phi$ takes $\hh^{\pm} \to \hh^{\pm}$ and ${\hh}^{\circ}$ to 
${\hh}^{\circ}$. As a consequence of this and the uniqueness of the factorization
(\ref{eqn123}),  $\phi$ fixes each of the three factors,
so $ x_{\pm} \in \hh^{\pm}_\phi$ and $ {x}^{\circ} \in {\hh}^{\circ}_\phi$.

Next we consider the map 
$ \Sigma: \hh^+_\phi \times \hh^{\circ}_\phi \times \hh^-_\phi \rightarrow  \mathcal{RB}\hh_\phi$.
We note that our argument above shows 
$\mathcal{RB}\hh_\phi = \hh_\phi \cap \mathcal{RB}\Lambda G$. 
Since $\hh_\phi$ carries the induced topology, $\mathcal{RB}\hh_\phi$ is open in $\hh_\phi$.
Moreover,
since the map $\Sigma$ is the restriction of the corresponding map for 
$\Lambda G$ and since $\hh_\phi$ is a Banach submanifold, it follows that $\Sigma$ 
 is an analytic diffeomorphism.
Finally, to prove that the image of the map under investigation is dense in $(\hh_{\phi})_I$,
we use the fact that, according to Theorem \ref{birkhoff},
 $\mathcal{RB}\Lambda G$ is given
as the complement of the zero set of a holomorphic section, $f_R$, on $(\Lambda G)_I$.
Hence $\mathcal{RB}\hh_\phi$ is the complement of the zero set of 
$f_R \big|_{\hh_\phi}$, and this 
section is real analytic. Therefore, $\mathcal{RB}\hh_\phi$ is either empty or dense,
and must be the latter, as it contains the identity.

\item
If $x$ is a non-constant loop in $\hh_\phi$, then, again, we have the equation (\ref{eqn123}).
Either $x_+$ or $x_-$ is a non-constant loop. Since $\phi$ is of the second kind, and
takes $\hh^\pm_1$ to $\hh^\mp_1$, it is impossible that such a non-constant factor be
fixed by $\phi$.  That is, the factors cannot both be in $\hh_\phi$.
\end{enumerate} \qed
 \end{proof}

\begin{remark} Many of the subgroups explicitly referred to in this 
article (such as the
$\sigma$-twisted subgroup $\Lambda G_\sigma$, where $\sigma$ is a complex linear
involution of the first kind, induced from an inner involution of $G$) are in fact isomorphic to $\Lambda G$.
However, it is important to note that the Birkhoff and Iwasawa decompositions we
prove for these are \emph{not} the same decompositions which would be obtained via this
isomorphism and the standard decompositions for $\Lambda G$. 
For example, with $\Lambda G_\sigma$,  the
three relevant subgroups, namely the constant subgroup and $\Lambda^\pm G_\sigma$,
do not correspond exactly to the corresponding subgroups $G$ and $\Lambda^\pm G$
under the isomorphism. 
\end{remark}

\subsection{$\tau$-Iwasawa factorization}  \label{sub6}
We saw that fixed point subgroups of automorphisms of the second kind are
not Birkhoff decomposable.  What we will later use in this
situation, for the case that $\tau$ is of order 2, is a more general
version of the Iwasawa factorization, which states that if $U$ is
a compact real form of $G$, then any element of $\Lambda G$ can be written
as a product $ug_+$, where $u \in \Lambda U$ and $g_+ \in \Lambda ^+G$.

It turns out that the discussion of (generalized) Iwasawa decomposable 
groups is much more complicated than the case of Birkhoff decomposable groups.
The reason for this will become clear, when we discuss examples below.

To keep the presentation as simple and short as possible we start from a very basic definition:
\begin{definition} \label{defniwasawadecomp1}
Let $\hh$ be some Banach Lie group, $\tau$ some real analytic involution and 
$\hh^\bullet$ some Banach Lie subgroup. Then $\hh$ together with the Banach Lie subgroups 
$\hh_{\tau}$ and $\hh^\bullet$ is called {\it $\tau$-Iwasawa decomposable} if  
$\hh_{\tau} \cdot \hh^\bullet$ is open in $\hh$.
\end{definition}
This definition even permits 
$\tau = id$ and $\hh^\bullet = \hh$.
However, we will apply this notion below exclusively in 
situations where $\tau$ is non-trivial, $\hh$ is a 
subgroup of $\Lambda G$ and $\hh^\bullet$ is $\hh^+$ (or, equivalently,
$\hh^-$).

\begin{remark}
1. As in the case of Birkhoff decomposable groups one obtains a representation
\begin{equation}
\hh = \bigcup_{\iota \in \mathcal{J}} \hh_{\tau} \cdot \iota \cdot \hh^\bullet
\end{equation}
where $ \mathcal{J}$ is a set of representatives for the orbits of the group action
 by 
$\hh_ {\tau} \times \hh^\bullet$ on $\hh$ given by $(h_ {\tau}, h_+ ) \cdot h = h_ {\tau} \cdot h h_+^{-1}$.
\end{remark}

2. Kellersch \cite{kellersch},
has investigated extensively Iwasawa decomposability in finite dimensional groups as 
well as in loop groups. He has given many examples illustrating what can happen.

3. If $\hh$ is a Birkhoff decomposable subgroup of some loop group $\Lambda G$ and 
$\tau$ an  involution of the second kind, then $\hh_{\tau} \cap \hh^\bullet \subset G$ and 
the question of obtaining a \emph{unique}
 $\tau$-Iwasawa decomposition reduces to 
finding a subgroup $\hh^\bullet$ such that this intersection is actually trivial.
This is a question for the finite dimensional group $G$, where it is known that
 one can in some cases
find complementary subgroups and in other cases one can not.

4. Kellersch \cite{kellersch} has given an example of a Birkhoff decomposable subgroup and an 
antilinear involution such that there are two disjoint open orbits. Thus, in general, one cannot assume that
$\hh_{\tau} \cdot \hh^\bullet$ is open and dense in $\hh$. In many examples this will be possible though.

To complete this round of definitions:
\begin{definition}
Let $\hh$, $\tau$ and $\hh^\bullet$ be as in Definition \ref{defniwasawadecomp1}.
An element $x \in \hh$ is called 
\emph{$\tau$-Iwasawa decomposable}
if it can be factorized 
\beq \label{decompeqn1}
x = z_\tau y_+,
\eeq
with $z_\tau \in \hh_{ \tau}$ and 
$y_+ \in  \hh^\bullet$.
\end{definition}

\begin{remark} \label{taudecompremark}
a) If $\hh$ is Iwasawa decomposable, then every element in some open neighbourhood of $I$ is
Iwasawa decomposable.

b) A $\tau$-Iwasawa decomposition 
$x = z_\tau y_+$ as above,
 can be written as 
$xy_+^{-1} = z_\tau 
= \tau z_\tau 
= \tau x  (\tau y_+)^{-1}$,
and the factorization (\ref{decompeqn1}) is equivalent
to the problem of finding $y_+ \in \hh ^\bullet$ such that
\beq  \label{requiredfact}
x^{-1} \tau x = y_+^{-1} \tau y_+.
\eeq
\end{remark}

\begin{theorem} \label{iwasawathm}
Suppose $\hh$ is any Birkhoff decomposable subgroup of an
arbitrary loop group $\Lambda G$,
and  $\tau$ is an involution of the second kind of $\hh$.
Set  $\hh^\bullet := \hh^+$. Then
\begin{enumerate}
\item
$\hh$ is  $\tau$-Iwasawa decomposable. 
\item Assume  $x ^{-1}\tau x  \in \mathcal{RB}\hh$ and Birkhoff decompose
 $ x ^{-1}\tau x  = v_+ \cdot a \ v_-$, where 
$v_{\pm} \in \hh^{\pm}_1$ and $ a \in {\hh}^{\circ}$. 
Then $x$ is $\tau$-Iwasawa decomposable if and only if $a= k^{-1} \, \tau k$ for some
$k \in {\hh}^{\circ}$.
\item
If the constant subgroup $\hh^\circ$ is \emph{connected} and \emph{compact} then 
every element, $x \in \hh$,
satisfying the requirement  that $x ^{-1}\tau x$ is in the identity 
component of  $\mathcal{RB}\hh$,
 is $\tau$-Iwasawa decomposable.
\item The left factor, $z_\tau$, in a $\tau$-Iwasawa decomposition 
(\ref{decompeqn1}) of $x \in \hh$, 
is unique up to right multiplication by a
constant  element of $\hh_{ \tau}$.
\item If $x$ is $\tau$-Iwasawa decomposable, then the factors $z_ \tau$ and $y_+$, 
in the decomposition $x = z_\tau y_+$,
 can, in a neighbourhood of $x$, be chosen to depend real analytically on $x$.
\end{enumerate}
\end{theorem}

\begin{proof}$(1)$ 
We need to show that $\hh _\tau \cdot \hh^+$ is open.
We apply example 3 of section 2.2.
Thus there exists an open  subset $\Omega$ of $\kk = \hh \times \hh$ such 
that every element of $\Omega$ can be written uniquely in the form
\beq \label{dbsplitting}
(h,f) = (u,u) \cdot (a,b) \cdot (v_+, v_-) ~ \in ~ \mathcal{K}_1^- \times \mathcal{K}^\circ \times \mathcal{K}_1^+. 
\eeq
Now we consider the extension of the given involution, $\tau (f,g) := ( \tau g, \tau f)$.
Then $\kk_{\tau} = \lbrace (h, \tau h), h \in \hh \rbrace $, and so we can identify $\hh$
with $\kk_\tau$.
If $x = (h, \tau h ) \in \kk_\tau$, then the two components of  (\ref{dbsplitting}) are
 $h=u a v_+$ and $\tau h = u b v_-$.  Applying $\tau$ to both  of these 
expressions, we also  have $\tau h = \tau u \, \tau a \, \tau v_+$ and
$h = \tau u \, \tau b \, \tau v_-$, to give the two decompositions
\beqas
(h, \tau h ) &=& (u,u) \cdot (a,b) \cdot (v_+, v_-) \\
&=& (\tau u, \tau u) \cdot (\tau b , \tau a) 
  \cdot (\tau v_- , \tau v_+) ~ \in  ~ \mathcal{K}_1^- \times \mathcal{K}^\circ \times \mathcal{K}_1^+.
\eeqas
Now the uniqueness of the factors
in this Birkhoff decomposition  implies that $\tau u = u$, and so
$\Omega \cap \kk_{\tau}$ is just the set $\{ (h, \tau h) = (u a v_+ , \tau(u a v_+)) | \tau u = u \} \subset \kk_\tau$. 
Since $\kk_ {\tau}$ carries the induced topology, $\Omega \cap \kk_{\tau}$ is open in $\kk_{\tau}$.
Finally, under our identification  of $\hh$ with $\kk_\tau$, this open set is nothing but the set
$\hh _\tau \cdot \hh^+$.

$(2)$ Since $ x ^{-1}\tau x \in \mathcal{RB}\hh$, we can write
 $ x ^{-1} \, \tau x = v_+ \, a \, v_-$, where 
$v_{\pm} \in \hh^{\pm}_1$ and $ a \in {\hh}^{\circ}$. 
Now the fact that $\tau(x^{-1} \, \tau x) = (x^{-1} \, \tau x)^{-1}$
 is then equivalent to
\bdm
v_-^{-1} \, a^{-1} \,  v_+^{-1} = \tau v_+ \, \tau a \, \tau v_-.
\edm 
By uniqueness of the factors in the Birkhoff decomposition, one has
 $v_-^{-1} = \tau v_+$ and $ \tau  a = a^{-1}.$
Hence
\beq \label{xinversetauxfact}
x^{-1} \, \tau x = v_+ \, a \, (\tau v_+)^{-1},
\eeq
and so the claim follows from Remark \ref{taudecompremark}.

$(3)$ Follows, in view of $(2)$, from  Proposition \ref{taufactprop} below. 

$(4)$ Assume we have two decompositions $zy_+ = hv_+$. Then 
\bdm
h^{-1}z = v_+ 
y_+^{-1}.
\edm
Clearly, the left element is fixed by $\tau$ and the right one is mapped 
to  $\hh^-$. Thus $ h^{-1}z$ is in the intersection of $\hh^+$ and $\hh^-$,
and is contained in $\hh_\tau \cap G$. This is the claim.

$(5)$ Also follows from Proposition \ref{taufactprop}, 
since the factors $v_+$ and $a$ in (\ref{xinversetauxfact}) both depend 
real analytically on $x$. \qed
 \end{proof}

To complete the proof of the Theorem above we show:
\begin{proposition} \label{taufactprop}
Let $H$ be a  semisimple Lie group and $\tau$ an  involution
of $H$. Define 
\beqas
H^* := \lbrace x \in H; ~\tau x = x^{-1} \rbrace, \\  
\mathcal{B} H^* := \{x \in H^* ; ~ \exists ~k \in H:~ x = k \, (\tau k)^{-1} \}.
\eeqas
\begin{enumerate}
\item  $\mathcal{B} H^*$ is open, and, locally,
the element $k$ can always be chosen to depend real analytically on $x$.
\item
If $H$ is compact then $\mathcal{B} H^*$ is both open and closed.
\end{enumerate}
\end{proposition}
\begin{proof}
Let $H_\tau$ denote the fixed point subgroup of $\tau$. Consider the map
$f: H/H_\tau \to H^*$, given by $f(kH_\tau) = k ({\tau k)}^{-1}$.
To prove both (1) and (2),
 we show that $f$ is a one to one immersion with open image. 
Then, for a family $x(t) \in \mathcal{B} H^*$,
there is a real analytic family $f^{-1}(x(t)) \in H/H_\tau$, and, on any
contractible set in $H/H_\tau$,
this lifts to a real analytic family $k(t)$.
 
Firstly, $a (\tau a)^{-1} = b (\tau b)^{-1} \Leftrightarrow b^{-1} a = \tau (b^{-1} a)
\Leftrightarrow b^{-1} a \in \hh_\tau \Leftrightarrow a \hh_\tau = b \hh_\tau$,
which shows that $f$ is one to one.

Secondly, $f$ is an
immersion, which can be seen as follows:
the tangent space to $H/H_{\tau}$ at $[k] := kH_\tau$ can be written in the form
$T_{[k]} (H/H_{\tau}) = \lbrace  r V ~;~ r \in \real, ~V \in \textup{Lie} (H), ~\tau V = -V \rbrace$, where we think of all groups as matrix groups and
$\tau$ is a linear map.
 To compute the differential of $f$ at $[k]$, 
 it is enough to  consider curves of the form
 $\gamma(t) = [k \exp (tV)]$, where $\tau V = -V$.
 Then $f(\gamma(t)) = k \exp (t V) \exp(tV)(\tau k)^{-1} = k \exp(2tV) (\tau k)^{-1} = k (\tau k)^{-1} + 2t k V (\tau k)^{-1} + o(t^2)$.
Hence $\dd f |_{[k]}(V) = 2 k V (\tau k)^{-1}$. The map $[V \mapsto 2 k V (\tau k)^{-1}]$ has zero kernel, so $\dd f |_{[k]}$  is injective.

Thirdly, the spaces have the same dimension,
because the tangent space at the identity of $H/H_\tau$ is equal to
the $-1$ eigenspace of $\tau$, which is the tangent space at the 
identity of $H^*$. Thus the image of $f$, namely $\mathcal{B} H^*$, is open. \qed
\end{proof}

\begin{remark} \label{doubleloopgroup}
Analysing the proof of the Theorem above one sees that an involution of the second kind
becomes an involution of the first kind, in the double loop group, for which we have 
already proven some properties.
Therefore, in some sense, only involutions of the first kind are of interest to us.
The double loop group set up is also appropriate for handling situations where 
two loops are necessary, rather than one.
\end{remark}


\section{Connection order $(_a^b)$ maps and their splittings}  \label{sec7}
In this section, we make the same assumptions on $G$ and the topology of $\Lambda G$
as in the previous sections.
Let $U$ be a manifold and define $\Lambda G(U)$ to be the set of smooth maps
$F: U \to \Lambda G$. Sometimes one wants to normalize $F$ to assume the 
value $I$ at some preassigned point $t_0 \in U$. In this case one says that $f$ 
is \emph{based (at $t_0$)}.
For any
extended integers $a$ and $b$ in ${\mathbb Z} \cup \{\pm \infty \}$,
with $a \leq b$, 
define  $\Lambda G(U)_a^b$ to be the subset of such maps
which are of connection order $(_a^b)$, that is,  whose Maurer-Cartan forms
have the Fourier expansion
\bdm
F^{-1} \dd F = A_a \lambda ^a + ... + A_b \lambda ^b,
\edm
for any $t \in U$.  
Define $\mathcal{B}\Lambda G(U)$ to be the subset of $\Lambda G(U)$ consisting of elements $F$
such that $F(t) \in \mathcal{B} \Lambda G$ for any $t$ in $U$, and $\mathcal{B}\Lambda G(U)_a^b$ 
analogously. 
Similar notation applies to subgroups and quotients of $\Lambda G$ in the obvious way.

\subsection{Splitting without an involution of the second kind} \label{sectswai}
\begin{proposition} \label{prop1}
 Let 
$\hh$ be any Birkhoff decomposable subgroup of $\Lambda G$, 
 and $a$ and
$b$ extended integers with $a<0 <b$. To every 
\bdm
F \in {\mathcal{B} \hh}(U)_a^b,
\edm
there corresponds a  unique pair
\beqas
G_- \in \hh(U)_a^{-1}, ~~~~~~ F_+ \in \hh(U)_1^b,
\eeqas
such that
\beqas
F= F_+F_- = G_-G_+
\eeqas
 for some $F_-$ and
 $G_+$ in $\hh ^-(U)$ and $\hh ^+(U)$ respectively.
If $F$ is assumed to be based at $t_0$,
 then  $F_+$, $F_-$, $G_+$, and  $G_-$ are also based at $t_0$.
\end{proposition}
\begin{proof}
$G_-$ and $F_+$ are just the left factors given by the left and right
Birkhoff decompositions respectively for $F$, with $F_+ \in \Lambda^+_1G$
and $G_- \in \Lambda^-_1 G$. 

 Let us check that 
$F_+ \in \hh(U)_1^b$:
For fixed $t \in U$, the Birkhoff factorization theorem states that
the map $(F_+,F_-) \mapsto F$ is a diffeomorphism.
 So we can deduce 
that the map $t \mapsto F(t) \mapsto F_+(t)$ is also differentiable.
Hence $F_+ \in \hh(U)$.
Now by definition, $F_+$ has the Fourier expansion
\beq \label{fplusexp}
F_+ = I + a_1 \lambda + a_2 \lambda ^2 + ...
\eeq
valid on $\overline{\DD}$, so 
\bdm
 F_+^{-1} \dd F_+ = b_1\lambda + b_2 \lambda^2 + ... 
\edm
is of bottom degree 1 in $\lambda$.
 For the top degree, note that $F_+$ also has the expression 
\beq \label{fplus}
F_+ = F F_-^{-1},
\eeq
and we are given that 
\bdm
F^{-1} \dd F = A_a\lambda^{a}+... + A_0 +...+ A_b \lambda^b.
\edm
Now comparing the Maurer-Cartan forms of both sides of (\ref{fplus}),
\bdm
 b_1\lambda + b_2 \lambda^2 + ... =
  F_-[A_{a}\lambda^{a} +...+ A_0 +...+ A_b \lambda^b] F_-^{-1}
  + F_- \, \dd (F_-^{-1}),
\edm
and using the fact that both $F_-$ and 
$F_-^{-1}$ are in $\Lambda ^- G$,  the conclusion is that both
 sides are  of top degree $b$ in $\lambda$.

The proof that $G_- \in \hh(U)_a^{-1}$ is similar.

Finally, assume that $F$ is based at $t_0 \in U$. Then 
$ I = F(t_0,\lambda)= F_+(t_0, \lambda) G_- (t_0, \lambda)$ for all $\lambda$.
Hence $F_+(t_0, \lambda) = F_0$ and $G_-(t_0, \lambda) = G_0$ are constant loops, and inverse to each other.
By the expression (\ref{fplusexp}) for $F_+$, it follows that $F_0 = G_0 = I$. \qed
\end{proof}

The following result is a little stronger than a converse 
to the last proposition:

\begin{proposition} \label{prop2}
Let $\hh$ be any Birkhoff decomposable subgroup of $\Lambda G$.
 Suppose $a$ and $b$  are 
extended integers with
$a<0<b$. To every pair
\beqas
G_- \in \hh(U)_a^\infty, ~~~~~~~~
F_+ \in \hh(U)_{-\infty}^b,
\eeqas
 such that $F_+^{-1}G_- \in \mathcal{B} \hh(U)$,
there corresponds a unique element 
\bdm
F \in \mathcal{B} \hh(U)_a^b,
\edm
which satisfies the equations
\beqa \label{defeqn1}
F &= &F_+F_-  \\
&=& G_-G_+, \label{defeqn2}
\eeqa
for some $F_-$ and $G_+$ in $\hh ^-(U)$ and $\hh ^+(U)$ respectively.
Moreover, if $F_+$ and $G_-$ are based at $t_0$, then $F$ is also based
at $t_0$.
\end{proposition}
\begin{proof}
By assumption, for each $t \in U$ there is a unique Birkhoff factorization
\beq \label{begineqn}
F_+^{-1}G_- = F_- G_+^{-1},
\eeq
where $F_-$ and $G_+$ are maps from $U$ to $ \hh ^-$ and $\hh ^+$ respectively and $F_-$ has constant term $I$.
  Rearranging this, we see that 
$F$ is well-defined by the equations (\ref{defeqn1}) and (\ref{defeqn2}).

Now calculate the Maurer-Cartan form for $F$, using (\ref{defeqn1}):
\bdm
F^{-1} \dd F = F_-^{-1}(F_+^{-1}\dd F_+)F_- + F_-^{-1}\dd F_-.
\edm
Because $F_+ \in \hh(U)_{-\infty}^b$ and $F_-$ has only powers
of $1/\lambda$,
this expression has top degree at most $b$ in $\lambda$, 
 and, using (\ref{defeqn2}), we similarly deduce that $F^{-1} \dd F$ also
has bottom degree $a$ in $\lambda$. In other words, $F \in \hh(U)_a^b$. 
Assume now that $F_+$ and $G_-$ are based at $t = t_0$, then
we have
\beqas
F(t_0, \lambda) = F_-(t_0, \lambda) = G_+(t_0, \lambda),
\eeqas
and this is a constant function of $\lambda$ by holomorphicity.
Since $F_-(t, \infty)= I$, this constant is the identity,
and so $F$ is also based at $t_0$. \qed
\end{proof}

\subsection{Splitting with an involution of the second kind} \label{sub8}
If $\hh$ is Birkhoff decomposable, then,
clearly,  Proposition \ref{prop1} applies to the case 
$F \in \mathcal{B} \hh_{ \tau}(U)_a^b$, where $\tau$ is an
involution of the second kind, although $G_-$ and $F_+$ will
not be fixed by $\tau$ themselves.  In fact, if $\tau F = F$,
we have
\bdm
\tau F_+ \, \tau  F_- = G_- \, G_+.
\edm
Since $\tau $ interchanges $\pm$, this presents 
two representations of an element in $\hh^-\cdot \hh^+$.
As the constant terms in the expansions of both $F_+$ and $G_-$
are taken to be the identity, uniqueness of the left Birkhoff factorization 
implies that $G_- = \tau F_+$.

With this in mind, we can state the analogue of Proposition \ref{prop2}
for a loop group with an involution, $\tau$, of the second kind. 
This depends on the $\tau$-Iwasawa factorization result,
Theorem \ref{iwasawathm},
and for this reason we define the \emph{Iwasawa big cells},
\bdm
\mathcal{B}^\tau_+ \hh = \hh_\tau \cdot \hh^+, \hspace{1cm} \mathcal{B}^\tau_-  \hh= \hh_\tau \cdot \hh^-,
\edm
and $\mathcal{B}^\tau_\pm \hh (U)_a^b$ to be the set of connection order $(_a^b)$
maps from $U$ which take values in $\mathcal{B}^\tau_\pm \hh$.
This is the exact analogue to the Birkhoff case, except that 
in general the Iwasawa cells are not dense. For many cases, they have
been investigated in \cite{kellersch}.
Further, if $\tau$ is an anti-linear involution defining a compact real
form, then one can use the standard
Iwasawa decomposition  to obtain the result without restrictions on $F$
globally, as in \cite{dorfmeisterpeditwu}.

\begin{remark} \label{switchremark}
The map $x(\lambda) \mapsto x(\lambda^{-1})$ is an analytic isomorphism of $\Lambda G$. It
follows that all the earlier definitions and results about Iwasawa decompositions 
are also valid if we switch
$+$ and $-$ and $\lambda=0$ with $\lambda =\infty$.
\end{remark}

\begin{proposition} \label{prop3}
Let  $\hh$ be any Birkhoff decomposable
subgroup of $\Lambda G$, and
 $\tau$ an  involution of the second kind of $\hh$.
Suppose $b$ is an  extended integer greater than zero.
Given an  element 
\bdm
F_+ \in \mathcal{B}^\tau_- \hh(U)_{-\infty}^b
\edm
there exists a unique element 
\bdm
[F] \in \frac{\hh_\tau}{\hh^\circ_\tau} (U)_{-b}^b
\edm
 such that for any local frame for $[F]$, that is, an element $F \in \hh_\tau(W)_{-b}^b$, 
 which represents $[F](t)$ for all $t$ in $W \subset U$,  one has the decomposition
\beqa \label{ftau1}
F &=& F_+ \,F_- \\
&=& \tau F_+ \, \tau F_-, \label{ftau2}
\eeqa
for some $F_-$  in $\hh ^-(W)$. Moreover, if $F_+$ is based at $t_0$, then
$F$ and $F_-$ can be chosen to be based at $t_0$ also.
\end{proposition}

\begin{proof}
By the definition of the Iwasawa big cell $\mathcal{B}^\tau_- \hh$,  one has, 
for each $t \in U$
the Iwasawa decomposition
\bdm 
F_+ = F  F_-^{-1},
\edm
where $F(t) \in \hh_\tau$. By Theorem \ref{iwasawathm} (combined with Remark \ref{switchremark}),
$F(t)$ is unique up to right multiplication by
$\hh^\circ_\tau$, which makes $[F](t)$ unique.
 By Item (5) of the same theorem, $F$ can locally be chosen as
a smooth function of $t$, and therefore the equivalence class
$[F] \in \hh_\tau/\hh^\circ_\tau (U)$ is smooth.
Since $F$ is fixed by $\tau$, we have $F = F_+ F_- = \tau F_+  \tau F_-$.
To check that $F$ is of connection order $(_{-b}^b)$, note that it follows from the 
expression $F=F_+ F_-$, where $F_+^{-1} \dd F_+$ is of top degree $b$ and 
$F_-$ takes values in $\hh^-$,
that $F^{-1}\dd F = ... + \alpha_b \lambda^b$ is of top degree $b$. Then
$F^{-1}\dd F = \tau (F^{-1}\dd F) = \beta_{-b} \lambda^{-b} + ...$ is of bottom degree $-b$,
due to the assumption that $\tau$ is an involution of the second kind. 
Finally, if $F_+(t_0) = I$, then $I=I \cdot I$ is a valid Iwasawa decomposition,
so one can choose $F(t_0) = F_- (t_0) = I$. \qed
\end{proof}

\subsection{Splitting with more than one  involution of the second kind} \label{sub9}
If we need to deal with the fixed point set of several (mutually commuting)
involutions of the second kind, 
then
we can always reduce this to the case of one involution of the second kind,
together with one or more additional involutions of the first kind.
It is straightforward to check the following facts:
\begin{lemma} \label{2lemma}
Let $\hh$ be any subgroup of $\Lambda G$,
 and suppose $\tau$ and $\hat{\tau}$ are a pair
of commuting involutions of the second kind of $\hh$. Define an 
involution
\bdm
\rho := \tau \circ \hat{\tau}.
\edm
Then 
\begin{enumerate}
\item
$\rho$ is an  involution of the first kind.
\item
If $\hat{\tau}$ is complex anti-linear and $\tau$ is complex linear
then $\rho$ is complex anti-linear.
\item
$\hh_{ \tau \hat{\tau}} = \hh_{ \rho \hat{\tau}} = \hh_{ \rho \tau}$.
\end{enumerate}
 \end{lemma}

\section{Potentials} \label{potentials}

In many applications of geometric interest one wants to use the method presented in this paper
to construct examples with specific properties.
In the context of the previous section one could say that one tries to find appropriate 
pairs of maps $F_+$ and $G_-$ as in Proposition \ref{prop2} or one map $F_+$ as in 
Proposition \ref{prop3}. The results just quoted then produce an ``extended frame $F$", from which the geometric object is derived  (more or less directly: see, for example, \cite{dorfmeisterpeditwu} and also Sections 6 and 7 below).
While, in general, the basic maps mentioned above are simpler than the frame $F$ constructed from them,
it is not trivial at all to write down those basic maps either. But one can simplify the situation by one more step:
one can consider the infinitesimal versions of the basic maps.
\begin{definition} 
 Let $\hh$ be any Birkhoff decomposable subgroup of $\Lambda G$, 
 and $a$ and $b$ extended integers with $a<0 <b$. Let $ F \in {\mathcal{B} \hh}(U)_a^b$ and
$G_- \in \hh(U)_a^{-1}$ and $ F_+ \in \hh(U)_1^b$,
such that $F= F_+F_- = G_-G_+$, and $F_+ |_{\lambda =0} = G_- |_{\lambda =\infty} = I$.
Then the pair  $(G_-, F_+)$ will be called a \emph{basic pair} for $F$,
 and the pair of Maurer-Cartan forms
\begin{equation}
\eta = ( G_-^{-1} dG_- , F_+^{-1} dF_+) = (\eta_- , \eta _+)
\end{equation}
will be called a {\it potential}
or  a {\it pair of potentials} for $F$. 
Similarly, for $\hh_\tau$, where $\tau$ is an involution of the second kind,
and $\eta_- = \tau (\eta_+)$,
as in Proposition \ref{prop3}, we call $F_+$ a \emph{basic map}, and its Maurer-Cartan
form $\eta_+$ a \emph{potential}.
\end{definition}

\begin{remark}

\begin{enumerate}

\item Of course, $\eta$ satisfies the obvious integrability condition. In some cases, such as 
\cite{dorfmeisterpeditwu}, the integrability for $\eta$ is trivial, while the integrability condition for
$F^{-1} dF$ is highly non-trivial.

\item In most cases one starts from some frame $F$ which is smooth 
on the manifold $U$.
The Birkhoff splitting used to derive the basic maps is not global and thus frequently 
introduces singularities.
In some cases, such as in \cite{dorfmeisterpeditwu}, one can modify the construction of the 
basic pairs so that one obtains smooth basic pairs. In general it is not known, whether for 
a smooth $F$ one can always find modified smooth basic pairs.

\item One should consider $\eta$ as a global object on $U$ which has singularities.

\item  What, in the case of  harmonic maps, has been called the DPW method is as follows:
start conversely, from some potential $\eta$ which has the proper expansion 
relative to $\lambda$ and satisfies the required integrability condition. Next one will 
integrate the equations
$\dd G_- = G_- \eta_-$ and $ \dd F_+ = F_+ \eta_+$, with some initial condition.
 If $\eta$ is free of singularities and if $U$ is simply connected, these differential 
equations have global solutions. The corresponding frame constructed in
Section \ref{sec7}, will,
 in general, still have singularities, but, wherever it is smooth,
  it will have the right properties.
At any rate, starting from some potentials with singularities one obtains geometric frames 
with singularities and vice versa.

\item In some geometric cases it is necessary that singularities occur,
 such as the case of constant negative Gauss curvature  surfaces in Euclidean
3-space, where Hilbert proved that there are no complete immersions.
For the case of $\hh_\tau$,
when the underlying Lie group is not compact, so that the Iwasawa decomposition
is not global, one can start with a potential which is smooth everywhere, and such
that the basic map $F_+$ is everywhere regular, but the extended frame $F$ has
singularities at points where the Iwasawa decomposition breaks down. This is investigated
for constant mean curvature surfaces in Minkowski space in \cite{brs}.

\item Clearly, the question of whether singularities can be removed 
or appear necessarily is a major problem that needs to be discussed.

\item Sometimes creating frames from potentials is confused with dressing 
(see the corresponding section below).
We would therefore like to point out that knowing all potentials is 
(on a theoretical level) equivalent with knowing all frames.
Dressing only creates new surfaces from old ones. It is a group action on the space of frames.
In all known cases, the dressing action  has many different orbits, and one usually does not 
know representatives for all these orbits. On the other hand  finding potentials which reflect,
to some extent, the properties one is looking for and dressing the corresponding frames has 
frequently lead to good results.
\end{enumerate}
\end{remark}



\section{Dressing solutions}  \label{dressing}
One of the methods for producing  new solutions of connection order
$(_a^b)$ maps  from given solutions  has been the dressing action of Zakharov and Shabat
\cite{zs}.  It was first discussed in terms of harmonic maps into 
Lie groups by Uhlenbeck  \cite{uhlenbeck1989}.
See \cite{guest} and \cite{terng2002} for more details and references.
In the case of Harmonic maps, a more general dressing action by subgroups
of $\Lambda \textup{GL}(n, {\mathbb C})$ is treated in \cite{bergveltguest}.

\subsection{Dressing for connection order $(_1^b)$ maps} \label{subsectdco}
If $\hh$ is Birkhoff decomposable, there is a natural local action of  
$\hh ^-$ on $\hh(U)_1^b$, where $U$ is some manifold, (and analogously of 
$\hh ^+$ on $\hh(U)_a^{-1}$) defined as follows: 
if $F_+ \in \hh(U)_1^b$  and $g_- \in \hh ^-$, then set 
$U^* = \lbrace t \in U ~|~ g_- F_+(t) \in \mathcal{RB}\hh \rbrace$.
Then,  for each $t \in U^*$, define a new
element  $\hat{F}_+(t)$ in $\hh_1 ^+$  by the
Birkhoff factorization
\bdm
g_-(\lambda) F_+(t,\lambda) = \hat{F}_+(t, \lambda) h_-(t, \lambda),
\edm

Applying Proposition \ref{prop1} to the map $F = g_- F_+$, which is
clearly of type $(_{-\infty}^b)$,  it follows that $\hat{F}_+$ is an
element of $\hh(U)_1^b$. It is straightforward to check that
this defines a left action by $\hh ^-$ on $\hh(U)_1^b$, or more
precisely a local action as one has to replace $U$ with $U^*$,
which is, in general, some large open subset of $U$.
Moreover, if $F_+$ is based at $t_0$, then  $\hat{F}_+$ is also based
at $t_0$. 

\subsection{Dressing for general order $(_a^b)$ maps} \label{sub10}
The above dressing action extends in a natural way to a left action
of $\hh ^- \times \hh ^+$ on $\hh(U)_a^b$, where $a<0<b$. If
$(g_-, g_+) \in \hh^- \times \hh^+ $ and $F \in \hh(U)_a^b$,
then $\hat{F} := (g_-, g_+)\diamond (F)$ is obtained on some subset 
$U^*$ of $U$ by first solving the
pair of equations
\beqas
g_-F = \hat{F}_+ h_-,\\
g_+F = \hat{G}_-h_+,
\eeqas
where the right hand sides are right and left Birkhoff factorizations
respectively.  If $F$ is of order $(_a^b)$ then $\hat{F}_+$ 
and $\hat{G}_-$ are of order
$(_1^b)$ and $(_a^{-1})$ respectively, so, applying Proposition
\ref{prop2} we can find a new element $\hat{F} \in \hh(U)_a^b$
which satisfies
\beqas
\hat{F} &=& \hat{F}_+ \hat{F}_-\\
 &=& \hat{G}_-\hat{G}_+,
\eeqas
for some $\hat{F}_-$ and $\hat{G}_+$.

Bearing in mind the principle that $F$ is associated to a pair $F_+ \in \hh(U)_1^b$
and $G_- \in \hh (U) _a^{-1}$ via Propositions \ref{prop1} and \ref{prop2},
it is unsurprising that this is the same action which one obtains by applying the 
action described above of $g_-$ to $F_+$ and $g_+$ to $G_-$ separately,
as the reader may easily verify. 

The action  can also be rephrased simply by using Example 3 of Section \ref{birkhoffsubsect}:
 what was described above 
simply is the dressing action  of $\kk^+$ on $(F,F) \in \kk^-_1$:
\begin{equation}
(h^+, h^-) (F,F) = (\hat{F},\hat{F})(v_+, v_-).
\end{equation}

\section{Isometric immersions of space forms into a sphere} \label{spaceforms}

Here we introduce the loop group formulation for the problem of 
isometric immersions of space forms into a sphere. The formulation
for immersions into hyperbolic space is analogous, replacing the 
group $SO(n+k+1)$ with the Lorentz group $SO_{-1}(n+k+1)$, and we 
describe that in Section \ref{hypsect} below.
In this section  we are primarily summarizing results from
\cite{feruspedit1996} and \cite{lms};
the key result is Proposition \ref{corlemma}.

\subsection{Extended frames and connections} \label{framesubsection}

As pointed out just above, we are interested in all immersions $f: M_c \rightarrow {\mathbb S}^n$,
where $M_c$ denotes some space form of curvature $c$ and 
${\mathbb S}^n$ denotes the unit sphere in $\real^{n+1}$. Wherever we say ``immersion", the map
is assumed to be of class $\mathcal{C}^\infty$.
Since such immersions do not exist in all cases globally, but on large submanifolds
which one obtains by removing ``singular points of $f$" from $M_c$,
we consider submanifolds $U \subset M_c$,
 where one can think of $U$ as being almost all of $M_c$.
For convenience we can even assume $U$ to be simply connected without changing 
the point of view above. In the discussion of adapted frames below, one considers
maps into a homogeneous space for which a global frame need not necessarily exist.
If one is interested in global problems, then
this can be handled in the present setting by considering maps into the the homogeneous
space $\hh/\hh^\circ$, as in \cite{brander2}. For convenience, we instead assume that
\emph{$U$ is a contractible open subset of $\real ^n$, with basepoint $t_0$}.

Let now $f: U \to {\mathbb S}^{n+k}$  be an immersion.  Then we will say $F: U \to SO(n+k+1)$ is an \emph{adapted frame}
for $f$ if the $(n+1)$'st column of $F$ is $f$ and if the derivative
$\dd f$ has no component in the directions given by the last $k+1$
columns.  This means that if $f$ is an immersion and 
\bdm
 F = [e_1,...,e_n,\xi _1,...,\xi _{k+1}],
\edm
(where $\xi_1 = f$), then $\{ e_i \}$ and
  $\{ \xi_i \}$ are orthonormal bases for the tangent and normal
spaces respectively of $f(U) \subset {\mathbb R }^{n+k+1}.$ Note that
an adapted frame always exists on a contractible set $U$.

In general, there is no canonical
choice of adapted frame for a given map $f$. 
The freedom is characterized as follows:
consider  the order two automorphisms of $SO(n+k+1,{\mathbb C})$ 
given by
\beqa
\sigma_2 = \textup{Ad}_P,  \label{sigma2} \\
\tau_2 = \textup{Ad}_Q, \label{tau2}
\eeqa
for the matrices 
\beqas
P := \bbar {cc} I_{n \times n} &0 \\ 0 & -I_{(k+1) \times (k+1)} \ebar, ~~~~~
Q := \bbar {cc} I_{(n +1) \times (n+1)} & 0 \\ 0 & -I_{k \times k} \ebar.
\eeqas
Here $I_r$ denotes an $r \times r$ identity matrix.
Let $SO(n+k+1,{\mathbb R})_{\sigma \tau}$ denote the elements of
$SO(n+k+1,{\mathbb R})$ fixed by both $\sigma_2$ and $\tau_2$.

\begin{lemma} \label{changeframe}
If $F: U \to SO(n+k+1,{\mathbb R})$ is an adapted frame for $f:U \to {\mathbb S}^{n+k}$,
then any other adapted frame for $f$ is given by
\beq \label{changeframeeq}
\hat{F} = FT,
\eeq
where $T: U \to SO(n+k+1,{\mathbb R})_{\sigma \tau}$.
\end{lemma}
\begin{proof}
A change of oriented orthonormal bases for the tangent and normal frames,
 keeping the $(n+1)$'st column (which is the map $f$) fixed,
amounts to right multiplication by a matrix of the form
\bdm
T = \bbar {ccc} 
  T_1 & 0 & 0 \\
  0 & 1 & 0 \\
  0 & 0 & T_2
\ebar ,
\edm
with $T_1 \in SO(n,{\mathbb R})$ and $T_2 \in SO(k,{\mathbb R})$.
But these are precisely the elements of  $SO(n+k+1,{\mathbb R})$ which
are fixed by both $\sigma_2$ and $\tau_2$. \qed
 \end{proof}

 Let $\Omega^1(U)$ be the space of real-valued one-forms on $U$
and denote by $A= F^{-1} \dd F \in so(n+k+1) \otimes \Omega^1(U)$ the 
pull-back by $F$ of the left 
Maurer-Cartan form of $SO(n+k+1)$. More explicitly,
\beqa  \label{adaptedform}
A  =    \bbar {cc}
      \omega & \betahat  \\
     -\betahat ^T  & \eta \\
  \ebar,
\eeqa
where $\omega \in so(n) \otimes \Omega^1(U)$ and
 $\eta \in so(k+1) \otimes \Omega^1(U)$.
If $f$ is an immersion, then $\omega$, $\eta$ and  $\betahat$ 
are the Levi-Civita connection form,
 the normal connection form and the second fundamental form
respectively of the immersion $f$, considered as a 
submanifold of  ${\mathbb{R}}^{n+k+1}$.
The requirement that $\dd f$ has no component in the direction
of any of the last $k+1$ columns
 is equivalent to 
the statement that the first row and first column of $\eta$ consist of
zeros. 

If $\mathfrak{g}$ is the Lie algebra of a Lie group $G$, and $A$ is
any 1-form in  $\mathfrak{g} \otimes \Omega^1(U)$ on the simply connected manifold $U$, then the 
condition for the (global) existence of a map $F$ into $G$ with 
connection 1-form $A$ is  the Maurer-Cartan equation 
\bdm 
\dd A + A \wedge A = 0.
\edm
Thus adapted frames are in one-to-one correspondence with 1-forms
$A \in so(n+k+1) \otimes \Omega^1(U)$ which satisfy the Maurer-Cartan
equation and have the property
that the first row and column of $\eta$ in (\ref{adaptedform})
are zero. We will call such one-forms \emph{adapted connection forms}.

If we write 
\bdm
\betahat = \bbar {cc}
  \theta & \beta
  \ebar,
\edm
where $\theta$ is the column vector $[\theta _i := e_i^T  \dd f]$
and $\beta$ is the second fundamental form of the map into ${\mathbb S}^{n+k}$,
then  $f$  being  immersive is equivalent with the 1-forms $\theta _i$ 
being linearly independent.  In a neighbourhood of such a point,  the induced sectional
 curvatures on the immersed submanifold 
are of  constant value $c$ if and only if 
\beq \label{curvatureeqn}
\dd \omega + \omega \wedge \omega = c \theta \wedge \theta ^ T.
\eeq

\begin{remark}
1.  For the case $c<1$,
 the smallest possible codimension for even the local
existence of such an immersion is $k=n-1$, and in this case it is known
\cite{moore1} that the normal bundle is flat:
\beq \label{flatnormal}
\dd \eta + \eta \wedge \eta = 0.
\eeq
Therefore, {\it in the rest of this paper, we will impose the  condition (\ref{flatnormal}) 
for all cases}, as it comes
naturally with the integrable system described below and thus 
allows us to apply integrable systems methods. 

2. Since we are going to discuss methods to produce adapted frames (without
necessarily satisfying the immersion condition), we define a slight
 generalization
of an isometric immersion with flat normal bundle 
of a space form into ${\mathbb S}^{n+k}$; generically, and in the appropriate
codimension, these will be isometric immersions.
\end{remark}

\begin{definition}  A \emph{constant curvature $c$ map with flat normal
bundle} into ${\mathbb S}^{n+k}$ is a map
$f: U \to {\mathbb S}^{n+k}$  corresponding to an
adapted frame as described above,
 and  where the
equations (\ref{curvatureeqn}) and (\ref{flatnormal}) hold.
\end{definition}

It was shown, in \cite{feruspedit1996} for the case $c\neq 0$ and
\cite{lms} for the case $c=0$, that one can insert an (appropriate complex or real)
auxiliary 
parameter $\lambda$ into the 1-form $A$, such that
the assumption  that the Maurer-Cartan equation holds for all $\lambda$ 
is equivalent
not only to the global  existence of $F$ on the simply connected manifold $U$ mentioned above, but also to the
equations (\ref{curvatureeqn}) and (\ref{flatnormal}), where $c$ will depend 
on $\lambda$ in general. The relevant facts are outlined below. 

\subsection{Extended frames for flat immersions} \label{flatsect}
\begin{lemma} \label{lemma1} 
Suppose that
\beq \label{typeA}
A^{\lambda} = \bbar {cc}
                \omega &  \lambda \betahat \\
                -\lambda \betahat ^T & \eta \\
	      \ebar
\eeq
is a family of complex matrix-valued one-forms, where $\omega$, $\eta$
 and $\betahat$
are independent of $\lambda$.  Then $A^\lambda$ satisfies the
Maurer-Cartan equation for all $\lambda$ if and only if 
  \begin{enumerate}
  \item
  $A^{\lambda_0}$ satisfies the Maurer-Cartan equation for some $\lambda_0 \in \cstar$,
  \item
  $\dd \omega + \omega \wedge \omega = 0$ and
 \item
  $\dd \eta + \eta \wedge \eta =0$.   
  \end{enumerate}
\end{lemma}
\begin{remark}
Lemma \ref{lemma1} 
is verified by writing out the Maurer-Cartan
equation for $A^{\lambda}$ and comparing coefficients of like powers
of $\lambda$, and the same comment applies to Lemma \ref{lemma2} below.
\end{remark}

\begin{definition} \label{extendA} 
A \emph{type $A$ extended connection} on $U$ is a family of 
$so(n+k+1,{\mathbb R})$-valued 1-forms $A^{\lambda}$ of the form (\ref{typeA}),
where 
\begin{enumerate}
\item $A^{\lambda}$ is real valued for real values of $\lambda$,
\item
$\omega$ is $n \times n$, 
\item
$\omega = 0 \hspace{2mm} \mbox{and} \hspace{2mm} \eta =0$,
\item
$A^{\lambda}$ satisfies the Maurer-Cartan 
equation for all $\lambda$.
\end{enumerate}
\end{definition}
Note that a type $A$ extended connection yields
 a family of adapted connection forms parameterized by $\lambda \in \rstar$.
These integrate to give adapted frames whose $(n+1)$'st columns are
constant curvature 0 maps with flat normal bundle, 
or \emph{flat maps}. Assuming the initial condition
\bdm
F(t_0, \lambda) = I,
\edm
at some base point $t_0 \in U$, the family of adapted frames $F(\lambda,t)$ corresponding to the 
extended connection $A^\lambda$ is unique.

Conversely, given a flat map $f$ into ${\mathbb S}^{n+k}$, 
one can choose an adapted frame $F$ for $f$ and obtains
the Maurer Cartan form $A = F^{-1} dF$, which is of the form
\beq \label{typeA'}
A = \bbar {cc}
                \omega &   \betahat \\
                - \betahat ^T & \eta \\
	      \ebar
\eeq

Introducing $\lambda$ as in (\ref{typeA}) above verifies that this connection 
form is integrable for all $\lambda$.

On the other hand, it is known that for a  given  flat immersion $f$ into ${\mathbb S}^{n+k}$, 
one can choose an adapted frame $F$ for $f$ such that both the tangent and the normal bundles
are \emph{parallel}, so that $\omega = 0 $ and $\eta =0$. 
Given an initial condition $F(t_0) = I$,
this parallel frame is unique.
Multiplying the corresponding connection form
by a parameter $\lambda$ gives a unique type $A$ extended connection.

It is worth noting here that the transition from any adapted frame for
a flat map to one with parallel bundles can be done independently of $\lambda$, and this argument 
applies to any connection order $(_0^b)$ map,  gauging it to an
order $(_1^b)$ map:
\begin{lemma} \label{gaugelemma}
If $F(\lambda,t)$ is a family of adapted frames
whose Maurer-Cartan forms satisfy all the conditions of Definition
\ref{extendA} with the exception of (3), 
then there is a unique map
$G: U \to SO(n+k+1)$ with $G(t_0) = I$ such that $F(\lambda,t)G(t)$ 
is parallel, that is, its Maurer-Cartan form satisfies all the conditions
of Definition \ref{extendA}.
\end{lemma}
\begin{proof}
If $A= A_0 + A_1 \lambda$ then the Maurer-Cartan equation for $A$ is
\bdm
\dd A_0 + A_0 \wedge A_0 + (\dd A_1 + A_0 \wedge A_1 + A_1 \wedge A_0) \lambda
+ A_1 \wedge A_1 \lambda^2 = 0.
\edm
Since this holds for all $\lambda$, the 1-form $A_0$ itself satisfies the 
Maurer-Cartan equation and there exists a map $H: U \to SO(n+k+1)$
such that 
\bdm
A_0 = H^{-1} \dd H.
\edm
Now set $G= H^{-1}$ and we have
\beqas	
(FG)^{-1} \dd (FG) &=& H(F^{-1} \dd F) H^{-1} + H\dd (H^{-1}) \\
&=& HA_1 H^{-1} \lambda.
\eeqas
Since $A_0$ is of the form 
$\bbar {cc}
                \omega & 0 \\
                0& \eta \\
	      \ebar$, $H$ has the form
$\bbar {cc}
                G_1 & 0 \\
                0& G_2\\
	      \ebar$, and therefore
\bdm
(FG)^{-1} \dd (FG) =
\bbar {cc}
                0 &  \lambda G_1 \betahat G_2^T \\
                -\lambda G_2 \betahat ^T G_1^T& 0 \\
	      \ebar 
\edm
which is a type $A$ extended connection. \qed
 \end{proof}

\subsection{Extended frames for immersions of non-flat space forms}  \label{sub4}
Here we summarize results from \cite{feruspedit1996}.
\begin{lemma} \label{lemma2} 
Suppose that 
\beq \label{typeB}
B^{\lambda} = \bbar {cc}
         \omega & \bbar {cc} (\lambda + \lambda^{-1}) \theta ~&~~ (\lambda - \lambda^{-1}) \beta \ebar \\
         -\bbar {cc} (\lambda + \lambda^{-1}) \theta ~&~~ (\lambda - \lambda^{-1}) \beta \ebar^T  & \eta 
	      \ebar
\eeq
is a family of complex matrix-valued one-forms, where $\omega$, $\theta$,
$\eta$  and $\beta$ are independent of $\lambda$.  Then $B^\lambda$ satisfies the
Maurer-Cartan equation for all $\lambda$ if and only if
 \begin{enumerate}
  \item
  $B^{\lambda_0}$ satisfies the Maurer-Cartan equation for some $\lambda_0 \in \cstar $, 
  \item    \label{curveqn}
$\dd \omega + \omega \wedge \omega =  4 \, \theta \wedge  \theta ^ T$.
  \item
  $\dd \eta + \eta \wedge \eta =0$.    \end{enumerate}
\end{lemma} 
\begin{remark}
In the case of isometric immersions, $\theta$ is an $n \times 1$ column
matrix and 
the coframe for the immersion $f$ will be 
\bdm
\hat{\theta} = (\lambda + \lambda ^{-1})\theta,
\edm
so  we will need that $\lambda \neq \pm i$ (for an immersion)
and  the condition 
at item \ref{curveqn} of the lemma means constant
sectional curvature 
\bdm c^\lambda = \frac{4}{(\lambda+\lambda^{-1})^2}. \edm 
\end{remark}

\begin{definition} \label{extendB} 
A \emph{type $B_j$ extended connection}
 on ${ \mathbb R}^n$ is a family of 
$so(n+k+1,{\mathbb C})$-valued 1-forms $B^{\lambda}$ of the form (\ref{typeB}),
where
\begin{enumerate}
\item
$\omega$ is $n \times n$,
\item
$\theta$ is $n \times 1$,
\item
the first row and column of $\eta$ are zero,
\item
 $B^{\lambda}$ satisfies the Maurer-Cartan equation for all $\lambda$ and
\item
the matrix coefficients of $B^{\lambda}$ satisfy the reality condition
 $\rbold _j$, chosen 
from the following three possibilities: 
\bdm
\begin{array} {cccc}
\rbold_1: &  (\rho_1 X)(\lambda) := \overline{X(-\bar{\lambda})}  = X(\lambda) 
 & \iff &  X(\lambda) \textup{real for }\lambda \in i \rstar,  \\
\rbold_2: &  (\rho_2 X)(\lambda) := \overline{X(\bar{\lambda})} = X(\lambda) 
 & \iff &  X(\lambda) \textup{real for} \lambda \in \rstar,  \\
\rbold_{-1}: &  (\rho_{-1} X)(\lambda) := \overline{X(1/\bar{\lambda})} = X(\lambda) 
& \iff & X(\lambda) \textup{real for} \lambda \in {\mathbb S}^1, 
\end{array}
\edm
\end{enumerate}
\end{definition}

\begin{remark}
As in the flat case, 
type $B_1$, $B_2$ and $B_{-1}$
extended connections correspond to
families of adapted connection forms parameterised this time by $\lambda$
 in  $i \rstar$,
$\rstar$ and ${\mathbb S}^1$ respectively.  Lemma \ref{lemma2} implies that the
corresponding maps $f^\lambda$ are constant curvature $c^\lambda$ maps,
with flat normal bundle, where
$c^\lambda = \frac{4}{(\lambda+\lambda^{-1})^2}$. Further, one has the following ranges for
$c^\lambda$:\\
\indent $\lambda \in i \rstar ~~~ \Rightarrow ~~~ c^\lambda \in (-\infty,0)$,\\
\indent $\lambda \in \rstar ~~~ \Rightarrow ~~~ c^\lambda \in (0,1]$,\\
\indent $\lambda \in {\mathbb S}^1 ~~~ \Rightarrow ~~~ c^\lambda \in (1, \infty) $.
\end{remark}
 \begin{definition} \emph{Type $A$} and \emph{type $B_j$ extended frames} are the
families of adapted frames obtained from  the type $A$ and type $B_j$
extended connections defined above, by integrating them with the initial
condition $F(t_0, \lambda) = I$ for each $\lambda$.
\end{definition}

Also as in the flat case, any constant curvature $c$ map, for $c \neq 0$,
 corresponds to
 a type $B_1$, $B_2$ or $B_{-1}$ extended frame 
depending on the value of $c$.
Unlike the flat case, where we could choose a parallel frame, here, even with the normalization 
$F(t_0, \lambda) = I$,
 the extended frame is only determined
up to a right multiplication by a matrix in $SO(n+k+1,{\mathbb R})_{\sigma \tau}$.

\subsection{Type $A$ extended frames with the reality condition $\rbold_1$} \label{sub5}
We defined extended frames of type $A$ as certain families of maps into
the Lie group with reality condition $\rbold_2$. 
However we could also have used the first reality condition
as follows: if
$A^\lambda = A_0 + A_1 \lambda$, with $A_i$ real and $\dd A + A \wedge A = 0$, 
define 
\beq  \label{correqn}
\hat{A}^{\lambda} := A_0 + i A_1 \lambda.
\eeq
 Then $\hat{A}$ also satisfies
the Maurer-Cartan equation and, for $\lambda = r_0 \in {\mathbb R}$, one has
\bdm
A^{r_0} = \hat{A}^{-ir_0}.
\edm
So $\hat{A}^{\lambda}$ integrates to a family of adapted frames of
flat maps for  $\lambda \in i \rstar$. We therefore define 
\emph{extended frames of type $A_1$ and $A_2$} accordingly.
 Note that extended frames of type $A$ with reality condition $\rbold_3$ do not exist.

Let $\{ e_k \}$ be the standard basis for ${\mathbb R}^{n+k+1}$.
If $I$ is one of the intervals 
\bdm
I_1 = (-\infty,0), \hspace{1cm} I_2 = (0,1), \hspace{1cm} I_{-1} = (1, \infty),
\edm
 let ${\mathbb S}_I(U,t_0)$ be the
set of families of constant curvature $c \in I$ maps with flat normal bundle
 $f: U \to {\mathbb S}^{n+k}$ such
that $f(t_0) = e_{n+1}$, and ${\mathbb S}_0(U,t_0)$ the analogous
set for flat maps.

We summarize the preceding discussion as:
\begin{proposition} \label{defprop}
\begin{enumerate}
\item
The elements of ${\mathbb S}_0(U,t_0)$ are in one to one correspondence
with type $A_j$ extended frames, where $j$ can be chosen to be either 1 or 2.
\item
The elements of ${\mathbb S}_{I_j}(U,t_0)$
are in one to one correspondence with equivalence classes of 
type $B_j$ extended frames, where $j$ is 1, 2, or -1.
The equivalence is up to right multiplication by a map 
$T: U \to SO(n+k+1,{\mathbb R})_{\sigma \tau}$.
\end{enumerate}
\end{proposition}

\begin{remark}
The statement in the second part  of the proposition above is frequently rephrased
as ``The elements of ${\mathbb S}_{I_j}(U,t_0)$
are in one to one correspondence with the gauge orbits of the $\lambda-$independent gauge group on the space of frames." 
\end{remark}

\subsection{Loop group formulation}  \label{loopformulation}
The above extended connections can be defined fairly naturally as 1-forms 
with values in subspaces of certain loop algebras as follows:
if $\mathfrak{g}$ is any Lie algebra, 
let $\widetilde{\mathfrak{g}} := \mathfrak{g} [\lambda, \lambda ^{-1}]$,
 the Lie algebra of Laurent polynomials with coefficients 
in $\mathfrak{g}$, which is given the Lie bracket 
\bdm  [\sum_{i=\alpha}^{\beta} X_i \lambda ^i,~
\sum_{i=\gamma}^{\delta} Y_i \lambda ^i] := \sum_{i,j} [X_i, Y_j] \lambda ^{i+j}.
\edm
 If $\widetilde{\mathfrak{H}}$ is any algebra of Laurent series 
in $\lambda$, then let us define $\widetilde{\mathfrak{H}}_a^b$ 
to be the vector subspace consisting of elements whose lowest power of
 $\lambda$ is $a$
and highest power of $\lambda$ is $b$.

Suppose we are given commuting 
Lie algebra automorphisms $\sigma_k$ and $\tau_2$ of 
$\mathfrak{g}$ of order $k$ and $2$ respectively, where $k \geq 2$.
 Using these,
define commuting automorphisms $\sigma$ and $\tau$ of the same order on 
$\widetilde{\mathfrak{g}}$ as follows:

\beqa \label{sigma}
(\sigma  X) (\lambda) := \sigma_k (X( (\zeta_k)^{-1}\lambda)), \\
(\tau  X) (\lambda) := \tau _2 (X(1/\lambda)), \label{tau}
\eeqa
where $\zeta_k$ is a primitive $k$'th root of unity.
Let $\gsigma$ be the subalgebra of $\widetilde{\mathfrak{g}}$ consisting of elements fixed by
$\sigma$, and define $\gsigmatau$ to be the subalgebra whose elements are fixed by both automorphisms. Assuming $\sigma$ and $\tau$ also commute with the involutions $\rho_i$,
we further define $\widetilde{\mathfrak{g}}_{\sigma}^{{\bf R}_j}$
 and $\widetilde{\mathfrak{g}}_{\sigma \tau}^{{\bf R}_j}$
to be the corresponding
 real subalgebras with the reality condition $\rbold_j$.

Now consider the case $\mathfrak{g} = so(n+k+1, {\mathbb C})$ and the 
inner automorphisms 
$\sigma_2 = \textup{Ad}_P$ and
$\tau_2 = \textup{Ad}_Q$
defined earlier by (\ref{sigma2}) and (\ref{tau2}).

It is straightforward to verify the following:
\begin{lemma} \label{loopalglemma}
\begin{enumerate}
\item  For $j=1,2$, type $A_j$ extended connections are precisely the one-forms 
which satisfy the Maurer-Cartan equation and whose
coefficients are in 
$(\widetilde{\mathfrak{g}}_{\sigma}^{{\bf R}_j})_1^1$.
 \item For $j = 1,2,3$ type $B_j$ extended connections are precisely the one-forms 
 which satisfy the Maurer-Cartan equation and whose 
 coefficients are in 
  $(\widetilde{\mathfrak{g}}_{\sigma \tau}^{{\bf R}_j})_{-1}^1$.
\end{enumerate}
\end{lemma}

Suppose  $\sigma_k$ and $\tau_2$ are commuting automorphisms of $G$, of
order $k$ and $2$ respectively,   and denote the
corresponding Lie algebra automorphisms with the same notation.
Extend $\sigma_k$ and $\tau_2$ to automorphisms $\sigma$ and 
$\tau$ of $\Lambda G$ by  (\ref{sigma}) and (\ref{tau}), and
define $\lgsigma$ and $\lgsigmatau$ to be the
 subgroups of $\Lambda G$ consisting of elements which are fixed, respectively,
  by $\sigma$, and
by both $\sigma$ and $\tau$.
Similarly we have the real subgroups $\lgsigmar$ and $\lgsigmataur$ where
${\rbold}$ is one of the reality conditions defined previously, with Lie
algebras $\Lambda \mathfrak{g}_\sigma ^\rbold$ and
 $\Lambda \mathfrak{g}_{\sigma \tau}^\rbold$ respectively. Thus 
Lemma \ref{loopalglemma} and Proposition \ref{defprop} can be rephrased as:
\begin{proposition} \label{corlemma}
\begin{enumerate}
\item
For $j = 1,~2$,
the elements of ${\mathbb S}_0(U,t_0)$ are in one to one
correspondence with the elements of $\Lambda G_{\sigma}^{\rbold_j}(U)_1^1$
that are based at $t_0$.
\item
Let $\hh = \Lambda G_{\sigma \tau} ^{\rbold_j}$.
For $j = 1,~2$ and $-1$, 
the elements of ${\mathbb S}_{I_j}(U,t_0)$ are in one to one
correspondence with the elements of $\frac{\hh}{\hh^\circ}(U)_{-1}^1$
that are based at $t_0$.
\end{enumerate}
\end{proposition}

\section{Correspondence between constant  positive,
negative and zero curvature  maps into a sphere or hyperbolic space}
\label{application}
We now  show how the above loop group splitting 
results can be used to  reduce the local problem of
constant curvature immersions with flat normal bundle into ${\mathbb S}^m$ and
${\mathbb H}^m$ to the flat case.

\subsection{Constant curvature immersions into hyperbolic space} 
\label{hypsect}

It turns out that the splitting results for connection order $(_a^b)$ maps applied
to immersions of space forms into the sphere lead naturally to 
immersions into hyperbolic space, in the case of the 
reality condition of the second kind, $\rbold_{-1}$.  We need the following discussion to
explain this. Define the Lie group
$SO_{-1} (n+k+1, {\mathbb C})$ to be the subgroup of $GL(n+k+1,{\mathbb C})$
consisting of matrices satisfying the equation
\beqas
X^T J X = J, \\
J := \textup{diag}(I_{n \times n}, -1, I_{k \times k}).
\eeqas
Extended frames for constant curvature
isometric immersions into hyperbolic space 
${\mathbb H}^{n+k}$ are defined analogously to those for the sphere 
 in Section \ref{flatsect}, replacing 
the Lie group $SO(n+k+1)$ with $SO_{-1}(n+k+1)$.  See \cite{feruspedit1996}
for a discussion of this in the non-flat case.
Define ${\mathbb H}_0(U,t_0)$ and ${\mathbb H}_I(U,t_0)$
analogously to their ${\mathbb S}$ counterparts,
 replacing the sphere with hyperbolic space.
The analogues of the results in Section \ref{spaceforms}  hold, 
replacing
${\mathbb S}$ with ${\mathbb H}$ and $G = SO(n+k+1,{\mathbb C})$ 
with $H = SO_{-1} (n+k+1, {\mathbb C})$; in addition, the formula for 
the induced curvature $c_\lambda$ is changed by a minus sign,
so that one must also replace $I_j$ with its reflection $-I_j$, for 
statements involving curvature in the interval $I_j$.  Thus the 
analogue to Proposition \ref{corlemma} is:

\begin{proposition} \label{corlemma2}
\begin{enumerate}
\item
For $j = 1,~2$, the elements of ${\mathbb H}_0(U,t_0)$ are in one to one
correspondence with the elements of $\Lambda H_{\sigma}^{\rbold_j}(U)_1^1$
that are based at $t_0$.
\item
Let $\hh = \Lambda H_{\sigma \tau}^{\rbold_j}$.
For $j = 1,~2$, or $-1$,
the elements of ${\mathbb H}_{-I_j}(U,t_0)$ are in one to one
correspondence with the elements of 
$\frac{\hh}{\hh^\circ}(U)_{-1}^1$ that are based at $t_0$.
\end{enumerate}
\end{proposition}

 Moreover, maps into $\Lambda G_\sigma$ with
a particular reality condition can be interpreted as maps into $\Lambda H_\sigma$
with a different reality condition via the next proposition. In order to state
the result more symmetrically, we introduce one further  reality
condition of the second kind:\\
\indent $\rbold_{-2}: ~~~~~  (\rho_{-2} X) (\lambda) := 
\overline{X(-1/\bar{\lambda})} =X(\lambda).$\\

\begin{proposition} \label{phiprop}
Consider the involutions of the first kind, defined by:
\beqa 
(\hat{\rho}_1 X )(\lambda) &:=& \tau \circ \rho_{-1} X(\lambda) \\
&=& Q \overline{X(\bar{\lambda})} Q^{-1}. \label{rhocon1}\\
(\hat{\rho}_2 X) (\lambda) &:=& \tau \circ \rho_{-2} X(\lambda) \\
&=&  Q \overline{X(-\bar{\lambda})} Q^{-1}, 
\eeqa
Define $\phi : GL(n+k+1, {\mathbb C}) \to GL(n+k+1, {\mathbb C})$ by 
\beqas
\phi X := \textup{Ad}_T X,\\
T := \textup{diag}(i I_{n \times n},~ 1, ~iI_{k \times k}).
\eeqas
Then 
\begin{enumerate}
\item
$\phi$ is an isomorphism between the subgroups
$G := SO(n+k+1,{\mathbb C})$ and $H:= SO_{-1} (n+k+1, {\mathbb C})$.
\item
$\phi$ defines a bijection
\bdm
 \Lambda G_{\sigma \hat{\rho}_j}(U)_1^1 ~\leftrightarrow~
  \Lambda H_{\sigma} ^{\rbold_j}(U)_1^1,
\edm for $j=1,2$.
\item
$\phi$ defines a bijection 
\bdm
\Lambda G_{\sigma \tau} ^{\rbold_j}(U)_{-1}^1 \leftrightarrow
 \Lambda H_{\sigma \tau} ^{\rbold_{-j}}(U)_{-1}^1,
\edm
 for $j=1,2$.
\end{enumerate}
\end{proposition}
\begin{proof} Let $\hat{F} := \phi F$.
\begin{enumerate}
\item Since $\textup{Ad}_T$ is a homomorphism, and clearly bijective,
it is enough to verify  that $F^T F = I$ is equivalent to
$\hat{F}^T J \hat {F} = J$ which is straightforward.
\item
Let $F$ be an element of $\Lambda G_{\sigma \hat{\rho}_j}(U)_1^1$.
To show that $\hat{F}$ is an element of
$\Lambda H_{\sigma} ^{\rbold_j}(U)_1^1$, we need to show
\begin{enumerate}
\item
$\hat{F}$ is of connection order $(_1^1)$,
\item
$(\rho_j \hat{F})(\lambda) = \hat{F}(\lambda),$
and
\item
$(\sigma \hat{F})(\lambda) = P\hat{F}(-\lambda)P^{-1} = \hat{F}(\lambda)$.
\end{enumerate}
All of these follow in a straightforward manner from the 
corresponding properties of $F$. We verify the reality conditions here
using the fact that  $F$ is fixed by both $\tau$ and $\hat{\rho}_j$, that is:
\beqas
F(\lambda )  = P F(-\lambda) P,\\
F(\lambda )  = Q\overline{F((-1)^{j+1} \bar{\lambda})}Q,
\eeqas
as well as  the readily verified identities
\beqas
PQT = -\bar{T}\\
 = -T^{-1}.
\eeqas
Thus we have
\beqas
(\rho_j \hat{F})(\lambda) &:=& \overline{\hat{F}((-1)^{j} \bar{\lambda)}}\\
&=& T^{-1} \overline{F((-1)^{j} \bar{\lambda})} T \\
&=& T^{-1} Q F(-\lambda)Q T \\
&=& T^{-1} Q P F(\lambda)PQT \\
&=& (-T) F(\lambda) (-T^{-1}) \\
&=& \hat{F}(\lambda)
\eeqas
\item
This is also straightforward, the main work being to check that if 
$F \in \Lambda G_{\sigma \tau} ^{\rbold_j}$ then $\hat{F} := T F T^{-1}$ is fixed
by both $\sigma$ and $\rho_{-j}$. 
\end{enumerate} \qed
 \end{proof}

Item (3) of Proposition \ref{phiprop} has the following immediate
corollary:
\begin{proposition} \label{simplecor}
 There are natural identifications
\beqas
 {\mathbb S}_{(-\infty,0)}(U,t_0) & \longleftrightarrow &
{\mathbb H}_{(-\infty,-1)}(U,t_0) \\
{\mathbb S}_{(1,\infty)}(U,t_0) & \longleftrightarrow &
{\mathbb H}_{(0,\infty)}(U,t_0) 
\eeqas
\end{proposition}

\subsection{Applying the splitting results to isometric immersions} \label{sub11}
We state the results locally, although similar statements would hold globally, away from
some singular set.
\begin{definition}
Let $k\geq n-1$.
The set of 
\emph{germs of constant curvature zero maps into ${\mathbb S}^{n+k}$,
at $t_0 \in U$} is defined by
\bdm
{\mathbb S}_0(t_0) := \Big( \bigcup _U {\mathbb S}_0(U,t_0) \Big) /\sim ,
\edm
where the union is over all neighbourhoods $U$ of $t_0$, and
two maps are equivalent if they agree on some neighbourhood of $t_0$.
\end{definition}

For $j= 1$, 2 or $-1$,
 let ${\mathbb S}_{I_j}(t_0)$ be the analogous set of germs
corresponding to families of constant curvature $c$ maps,
with $c \in I_j$.

\begin{theorem}  \label{cortheorem}
\begin{enumerate}
\item
 There are canonical bijections
\bdm
{\mathbb S}_{(-\infty,0)}(t_0)  \longleftrightarrow {\mathbb S}_{0}(t_0) \longleftrightarrow
{\mathbb S}_{(0, 1)}(t_0).
\edm
\item
 There is a
canonical bijection
\bdm
{\mathbb S}_{(1,\infty)}(t_0) \longleftrightarrow
 {\mathbb H}_{0} (t_0).
\edm
\end{enumerate}
\end{theorem}
\begin{proof}
\begin{enumerate}
\item
An element of ${\mathbb S}_{(-\infty,0)}(t_0)$ or
 ${\mathbb S}_{(0, 1)}(t_0)$
corresponds by Proposition \ref{corlemma} to an element 
$[F] \in \frac{\lgsigmataur}{(\lgsigmataur)^\circ}(W)_{-1}^1$,
based at $t_0$,
and where $W$ is some neighbourhood of $t_0$, and
 $\rbold$ is $\rbold_1$ or $\rbold_2$ respectively. By Proposition
\ref{prop1}, this corresponds to a unique element 
$F_+ \in \Lambda G_{\sigma}^{\rbold}(W)_1^1$, also based at $t_0$,
and which, according to
Proposition \ref{corlemma},  gives an element of
${\mathbb S}_0(t_0)$.

Conversely, given a based element $F_+ \in \Lambda G_{\sigma}\rbold(W)_1^1$,
using Proposition \ref{corlemma} again
and Proposition \ref{prop3}, there is a unique  element of
 $\frac{\lgsigmataur}{(\lgsigmataur)^\circ}(W)_{-1}^1$, based at $t_0$.
\item
For the case $c \in (1,\infty)$,
recall that constant curvature $c$ maps  are represented
by maps into the loop group with the reality condition
$\rbold_{-1}$, which is of the second kind and commutes with $\tau$.

Using Lemma \ref{2lemma}, the analogue of 
 part (1) of this
theorem  gives  local correspondences between
these constant curvature maps and based elements of 
$\Lambda G_{\sigma \hat{\rho}_1}(U)_1^1$, where $\hat{\rho}_1$ is the 
anti-linear involution of the \emph{first} kind  given by
(\ref{rhocon1}).\\
Using the bijection at item (2) of Proposition \ref{phiprop},
\bdm  
\Lambda G_{\sigma \hat{\rho}_1}(U)_1^1 ~\leftrightarrow~
\Lambda H_{\sigma} ^{\rbold_1}(U)_1^1,
\edm 
we obtain the required correspondence.
\end{enumerate} \qed
\end{proof}

Finally, the analogue of Theorem \ref{cortheorem} for immersions into
hyperbolic space, whose proof is essentially the same is:
\begin{theorem}  \label{cortheorem2}
\begin{enumerate}
\item
 There are canonical bijections
\bdm
{\mathbb H}_{(0,\infty)}(t_0)  \longleftrightarrow {\mathbb H}_{0}(t_0) \longleftrightarrow
{\mathbb H}_{(-1,0)}(t_0).
\edm
\item
 There is a
canonical bijection
\bdm
{\mathbb H}_{(-\infty,-1)}(t_0) \longleftrightarrow
 {\mathbb S}_{0} (t_0).
\edm
\end{enumerate}
\end{theorem}

\subsection{Example}
It is difficult to exhibit examples with explicit formulae: given 
an isometric immersion, there is a recipe for inserting the parameter
$\lambda$ into the Maurer-Cartan form of an adapted frame, but then
one can generally only integrate this Maurer-Cartan form numerically.
Further, the Birkhoff and Iwasawa splittings must generally be computed
numerically.  In the case of constant mean curvature
 surfaces, software has been written
which can produce images of the surface corresponding to any ``potential"
used in the DPW method.

Here we consider one simple example, from \cite{brander2},
 of a family of  2-dimensional spheres,
of constant curvature $c_\lambda \in (1,\infty)$, isometrically embedded 
in ${\mathbb S} ^3$,  and we compute the 
corresponding flat immersion. Consider the family of maps 
$F_\lambda : \real ^2 \to G = SO(4,\cc)$ which takes $(u,v) \in \real^2$
to the matrix
\bdm
 \bbar {cccc}
     \cos u ~&~ -\sin u \sin v ~&~ a \sin u  \cos v  ~&~ b \sin u  \cos v  \\
      0  & \cos v                & a \sin v          & b \sin v   \\
 -a \sin u  ~~&~~ -a\cos u \sin v  ~~&~~ a^2 \cos u \cos v  + b^2 ~~&~~ 
               ab (\cos u \cos v -1)\\
     -b\sin u  ~~&~~  -b\cos u \sin v  ~~&~~ ab(\cos u \cos v -1) ~~&~~
                 b^2\cos u \cos v +a^2   \ebar,
\edm
where 
\bdm
a= \frac{1}{2}(\lambda + \lambda^{-1}),  \hspace{1cm} 
b= \frac{i}{2}(\lambda - \lambda^{-1}).
\edm

The Maurer-Cartan form of $F_\lambda$ is
\bdm
F_\lambda ^{-1}\dd F_\lambda = \bbar {cccc}
   0 ~~&~~ -\sin v  \, \dd u ~~&~~ a \cos v  \, \dd u ~~&~~ b \cos v  \, \dd u \\
  \sin v  \, \dd u & 0 & a \, \dd v &   b \, \dd v \\
  - a \cos v  \, \dd u & -a \, \dd v & 0 & 0 \\
  -b \cos v  \,  \dd u & -b \, \dd v & 0 & 0
  \ebar.       
\edm
It is easy to check that $F_\lambda^{-1}\dd F_\lambda$
 is fixed by $\sigma$, $\tau$ and $\rho_{-1}$, so it
takes values in the Lie algebra of 
$\hh:= \Lambda G_{\sigma \tau}^{\rbold_{-1}}$.
Because $F_\lambda (0,0) = I \in \hh$, it 
follows that  $F_\lambda$ is a map into $\hh$,
and, since its Maurer-Cartan form has top and bottom degree 1 and -1
respectively, it represents an element of 
$\frac{\hh}{\hh^\circ}(\real^2)_{-1}^1$. 
Thus, according to Section \ref{spaceforms}, if the the third column $f^\lambda$
of $F_\lambda$
 is an immersion, then, for a value of $\lambda$ in ${\mathbb S}^1$, it is an immersion
into ${\mathbb S}^3$ with constant Gauss curvature greater or equal to 1.
The coframe for $f^\lambda$ is given by 
\beq \label{theta}
\theta = \bbar {c} \frac{1}{2}(\lambda + \lambda^{-1}) \cos v \, \dd u \\
           \frac{1}{2}(\lambda + \lambda^{-1}) \, \dd v \ebar,
\eeq
and so, if $\lambda \neq \pm i$, then $f^\lambda$ is
immersive away from the 
 degenerate coordinate lines $\cos v =0$.
In fact $f^\lambda$ is a deformation, through a family of isometrically
embedded spheres, of the totally geodesic embedding of
${\mathbb S}^2$ into ${\mathbb S}^3$ given by 
$f(u,v) = \bbar {cccc}  \sin  u \cos v, ~&~
          \sin v, ~&~
        \cos u \cos v, ~&~
         0 \ebar^T$, which is achieved at $\lambda = 1$.

To understand the corresponding flat immersion $f_+$ into ${\mathbb H}^3$,
given by the second part of Theorem \ref{cortheorem}, let 
\bdm
F = F_+ F_-
\edm
represent the right Birkhoff splitting of the frame $F$,
valid on some open set $U \subset \real^2$.  Then $F_+$ is an
adapted frame for $f_+$.  Expressing $F_+^{-1} \dd F_+$ in terms of $F^{-1} \dd F$ 
and the function $F_-$, we deduce that
\bdm
F_+^{-1} \dd F_+ = K \alpha_{+} K^{-1} \lambda,
\edm
where $F^{-1} \dd F = \alpha_- \lambda^{-1} + \alpha_0 + \alpha_+ \lambda$
and $K$ is a map into $\hh^\circ = SO(2) \times SO(2)$, represented by diagonal
matrices whose components are $2 \times 2$ blocks.

Thus we have 
\beqas
\alpha_+ = \bbar {cc} 0 ~&~ A \\ -A^T ~&~ 0 \ebar, \hspace{1cm}
A = \frac{1}{2} \bbar{cc} \cos  v \, \dd u ~& ~i \cos  v \, \dd u \\ \, \dd v ~&~ i \, \dd v \ebar,\\
K = \bbar{cccc} \cos \theta ~~&~~ \sin \theta ~~&~~ 0 ~~&~~ 0 \\ -\sin \theta ~~&~~ \cos \theta ~~&~~ 0 ~~&~~ 0 \\
    0  ~~&~~ 0  ~~&~~ \cos \phi  ~~&~~ \sin \phi \\ 0  ~~&~~ 0  ~~&~~ -\sin \phi  ~~&~~ \cos \phi \ebar ,
\eeqas
where  $\theta$ and $\phi$ are functions $U \to \real$.
Putting these together, we have 
\beqas
F_+^{-1} \dd F_+ = \bbar{cc} 0 ~&~ B \\ -B^T ~&~ 0 \ebar \lambda ~=:~ \eta_+ \lambda, \\
B = \frac{1}{2} \bbar{cc} e^{i\phi}( \cos \theta \cos v \, \dd u + \sin \theta \, \dd v)
  & i e^{i \phi} ( \cos \theta \cos v \, \dd u + \sin \theta \, \dd v) \\
  e^{i \phi} (-\sin \theta \cos v \, \dd u + \cos \theta \, \dd v ) ~~&~~
  i e^{i \phi} (-\sin \theta \cos v \, \dd u + \cos \theta \, \dd v ) \ebar.
\eeqas
Now $F_+^{-1} \dd F_+$ must  satisfy the Maurer-Cartan equation for all values of
$\lambda$, which is equivalent, in this case, to the pair of equations 
$\eta_+ \wedge \eta_+ = 0$ and $\dd \eta_+ = 0$. In particular, all the components
of the matrix $B$ can be integrated, and we can write
\bdm
B = \bbar{cc} \dd x ~&~ i \,  \dd x \\ \dd y & i \, \dd y \ebar,
\edm
for a pair of functions $x$ and $y$ which map $U \to \real$ (we take
$U$ to be simply connected).
It follows that, for a fixed value of $\lambda$ in $i \real$,
 $F_+$ is a part of the torus in $SO(4, \mathbb C)$ given by
exponentiating the Abelian subalgebra, $\mathfrak{m}$ of $so(4,\mathbb C)$,
given by the span over $\real$ of the two matrices 
\bdm
X = \bbar{cccc} 0  ~&~ 0  ~&~ 1  ~&~ i \\
                0  ~&~ 0  ~&~ 0  ~&~ 0\\
                -1  ~&~ 0  ~&~ 0  ~&~ 0\\
                -i  ~&~ 0  ~&~ 0  ~&~ 0\ebar  ,
                \hspace{1cm}
Y = \bbar{cccc} 0  ~&~ 0  ~&~ 0  ~&~ 0 \\
                0  ~&~ 0  ~&~ 1  ~&~ i\\
                0  ~&~ -1  ~&~ 0  ~&~ 0\\
                0  ~&~ -i  ~&~ 0  ~&~ 0\ebar  .
\edm
The surface in ${\mathbb H}^3$ which corresponds to the third column of $F_+$,
is the analogue in  ${\mathbb H}^3$ of the Clifford torus in ${\mathbb S}^3$. Explicitly,
$F_+ = \exp(x X \lambda) \exp(y Y \lambda))$, and the third column  of this
matrix is $[\, x \lambda, ~\,  y \lambda, ~\, \frac{1}{2}(2-x^2 \lambda ^2 - y^2 \lambda^2), ~\,
  -\frac{1}{2}i (x^2 \lambda^2 + y^2 \lambda^2) \, ]^T$.  After applying $\textup{Ad}_T$,
  as given in the proof of Proposition \ref{phiprop}, this becomes:
\bdm
f_+ = [\, i x \lambda, ~\, i y \lambda, ~\, \frac{1}{2}(2-x^2 \lambda ^2 - y^2 \lambda^2), ~\, 
  \frac{1}{2} (x^2 \lambda^2 + y^2 \lambda^2) \, ]^T. 
\edm
which indeed takes values in ${\mathbb H}^3$ for $\lambda \in i \real$.       

\begin{remark}
In the application of the DPW method to harmonic maps from a Riemann surface
to a symmetric space \cite{dorfmeisterpeditwu}, 
the basic map (corresponding to $F_+$ here) is a holomorphic map into the
complex loop group $\Lambda G$. The same construction works for
 pluriharmonic maps \cite{dorfech}.
 Note that in our application to space forms, $F_+$ maps into
some real form of $\Lambda G$, and is of course not holomorphic.
 However, as shown in \cite{brander2},
after complexifying the source manifold, $F_+$ can easily be extended
to a holomorphic map into  $\Lambda G$, provided $F_+$ is real analytic.
Taking this as the basic map in the reverse direction of the DPW method,
one then obtains an extended frame for a pluriharmonic map. For the example
exhibited here, this pluriharmonic map would be from some open subset in $\cc ^2$ 
into the symmetric space $\frac{SO(4)}{SO(2) \times SO(2)}$. 
\end{remark}

\providecommand{\bysame}{\leavevmode\hbox to3em{\hrulefill}\thinspace}
\providecommand{\MR}{\relax\ifhmode\unskip\space\fi MR }
\providecommand{\MRhref}[2]{%
  \href{http://www.ams.org/mathscinet-getitem?mr=#1}{#2}
}
\providecommand{\href}[2]{#2}

\end{document}